\documentclass[a4paper]{article}
\usepackage{RR}
\usepackage{hyperref}
\RRdate{Avril 2008}

\RRauthor{
H\'el\`ene Barucq%
\thanks[sfn]{EPI Magique-3D, Centre de Recherche
  Inria Bordeaux Sud-Ouest}
\thanks[sf1]{Laboratoire de Math\'ematiques et de
  leurs Applications, CNRS UMR-5142, Universit\'e de Pau et des Pays de
  l'Adour --  B\^atiment IPRA, avenue de
  l'Universit\'e  -- BP 1155-64013 PAU CEDEX}%
  \and
Bertrand Duquet\thanks{Institut Fran\c cais du
  P\'etrole, BP311, 92506 Rueil Malmaison Cedex,
  France}%
\and
Frank Prat\thanksref{sfn}\thanksref{sf1}
}
\authorhead{Barucq, Duquet \& Prat}
\RRtitle{Propagation one-way \`a amplitude pr\'serv\'ee
dans des milieux h\'et\'erog\`enes}
\RRetitle{True amplitude one-way propagation in heterogeneous media}
\titlehead{True amplitude one-way propagation}
\RRresume{Dans ce travail, on s'int鲥sse \`a deux formulations one-way de l'\'equation des ondes acoustique qui ont \'et\'e construites \`a partir du mod\`ele complet one-way propos\'e par M.V. De Hoop. Le principal objectif de cette \'etude est de comparer les deux formulations dans lesquelles on a introduit un terme permettant d'am\'eliorer le calcul de l'amplitude de la solution num\'erique. Il ressort de l'analyse que m\^eme si les deux syst\`emes sont \'equivalents du point de vue th\'eorique, il n'en est rien au niveau des performances num\'eriques. On montre en particulier que tant que le milieu de propagation est r\'egulier, les deux mod\`eles se comportent identiquement et que les diff\'erences sont nettes si le milieu comporte des h\'et\'erog\'eit\'es. On peut donc en conclure que le pr\'ecision de l'amplitude est tr\`es sensible \`a la formulation du mod\`ele.}
\RRabstract{This paper deals with the numerical analysis of two one-way systems derived from the general complete modeling proposed by M.V. De Hoop. The main goal of this work is to compare two different formulations in which a correcting term allows to improve the amplitude of the numerical solution. It comes out that even if the two systems are equivalent from a theoretical point of view, nothing of the kind is as far as the numerical simulation is concerned. Herein a numerical analysis is performed to show that as long as the propagation medium is smooth, both the models are equivalent but it is no more the case when the medium is associated to a quite strongly discontinuous velocity.}
\RRmotcle{Analyse micro-locale num\'erique, formulation one-way, ondes acoustiques}
\RRkeyword{Microlocal numerical analysis, one-way formulation, acoustic waves}
\RRprojets{Magique-3D}
\RRtheme{\THNum} 
\RCBordeaux

\RequirePackage{calc}
\usepackage{graphics,epsfig,graphicx,subfigure,pifont}
\usepackage{setspace,amsmath,amsfonts,amssymb,amsthm,delarray}
\usepackage{fancyhdr}
\usepackage{makeidx}   

\usepackage{pstcol}
\usepackage{pstricks}
\usepackage{fancybox}
\usepackage{enumerate}

\newtheorem{theoreme}{Theorem}[section]

\newtheorem{proposition}[theoreme]{Proposition}

\newtheorem{corollaire}[theoreme]{Corollary}

\newtheorem{remarqueit}{\it Remark}

\newcommand{\rr}{I\hspace{-0.3em}R}
\newcommand{\nn}{I\hspace{-0.3em}N}

\newcommand{\mcaption}[1]{\caption[]{\textit{#1}}}

\def\vec#1{\mathchoice{\mbox{\boldmath$\displaystyle#1$}}
{\mbox{\boldmath$\textstyle#1$}}
{\mbox{\boldmath$\scriptstyle#1$}}
{\mbox{\boldmath$\scriptscriptstyle#1$}}}

\def\directorywork{./figure/}

\begin{document}
\makeRR

\section{Introduction} The numerical solution of seismic acoustic waves
propagation in heterogeneous media is generally based on the solution of the second-order full wave
equation. Then the second-order wave equation can be solved
completely by using a finite difference scheme and it is
well-known that such an approach results in a high computational burden, especially
in the case of three-dimensional problems. This surely explains that a
lot of people prefer to solve an approximate problem which
involves either a truncate expansion of the solution \cite{Farra}
or an approximate wave equation arising from the factoring of the
exact one \cite{zhang,kiya}.  In the simplest case of a
homogeneous medium, the solution can be obtained from the
inversion of a system of one-way wave equations
\cite{aki}. It is a very interesting way of solving the
wave equation because it is based on the decomposition of the wave
into a down-going part and an up-going one which reproduces the
physical phenomenon very faithfully. During several years, people
(see \cite{Corones} and its references) tried to extend this
approach by introducing correcting terms in the model to account
for heterogeneities into the propagation medium. In 1996, M.V. De
Hoop \cite{ref1} derived a new formulation based on the
micro-local analysis which allowed to derive a complete system of
one-way wave equations coupling by exact correcting terms. The use
of the theory of pseudo-differential operators makes the
derivation of the system easy. Nevertheless, its plain writing is complicated
because it involves the composition of pseudo-differential
operators which means that each term is defined from an asymptotic
expansion. Thus in practice the complete system is approximated by truncating the asymptotic expansions. Such an approach may look as if it complicates the solution as compared to the now well-controlled solution of the full wave equation. But its formulation allows to unpack the multiples from the primary reflections which is of outstanding importance for the geological interpretations and allows to reduce the computation time. Paper \cite{ref1} has been followed by numerous
publications and among them, we refer to as \cite{prat} in
which one can find a complete bibliography on the topic. As far
as the numerical solution is concerned, it is associated to the
inversion of an approximate system generally based on only keeping the main term in each asymptotic development. J. Le Rousseau was the first to obtain accurate
snapshots (see \cite{ref3} and \cite{ref2}). However if using
his numerical method for the computation of arrival times, one
gets erroneous results on the amplitudes level. More recently, Zhang \emph{et al.} \cite{zhang} have proposed a corrected one-way wave equation which allows to compute the correct amplitude of the acoustic pressure. This new one-way wave equation is obtained from the factorization of the full wave problem. The equation is factorized by using a WKBJ solution in which the amplitude is taken into account as well as its phase. Their new formulation of the one-way system includes a correction term.\\
In this work, we
intend to show that the amplitude of the numerical solution can be
corrected by adding a transmission term in the system, proposed in \cite{ref3}. The heterogeneities of the medium are modeled as discontinuities of the velocity which is supposed to vary in all the directions. The correcting term
can be included in the system by two ways and we show that one
of them is optimal. In fact we show that the best improvement of the amplitude is obtained when including the transmission operator into the right-hand-side of the system. This might come upon the reader since the model which is the nearest of the one of Zhang \emph{and al.} \cite{zhang} is obtained by including the transmission operator into the one-way equations, \emph{i.e.} into the left-hand-side of the system. But some numerical tests indicate that the two approaches are close in the case of smooth media. \\
The paper is divided into 7 sections, plus this one and a concluding part. The next one 
deals with the initial model whose unknown is the acoustic
pressure and its transformation into a first-order system by
introducing the velocity as an unknown also. The third part is
devoted to the reduced system whose derivation is based on
selecting the depth variable as the leading direction and by plugging
the other terms into the frequencies domain after
 using a Fourier-Laplace transform. The fourth part concerns the first-order approximation of the reduced system. By accounting
 for the complete coupling terms of order 0, we get two equivalent systems which are interesting to consider since their
 numerical solution can be obtained by two different ways.
The fourth next parts deal with some numerical aspects which are essential for the method. We have chosen to neglect the description of the propagation because we intend to focus on the transmission operator. The numerical tests are developed in the 2D case but we mention that some 3D test have been performed in \cite{prat}.\\

In the following, we use standard notations for the micro-local analysis of classical pseudo-differential operators and we refer to \cite{T} for their definitions. We only precise that the symbol of an operator $P$ is denoted by $\sigma(P)$ and its principal symbol is $\sigma_P(P)$.

\section{Initial model}

The analysis of the waves propagation is an efficient tool for imaging the soil. Assume the region of interest $\Omega$ is located between the surface of the earth $\{z=0\}$ and a given depth $\{z=z_{max}\}$. The phenomenon of propagation is supposed to be generated by a source located at the top of $\Omega$. Then the discontinuities of the propagation velocity can be defined by computing the time arrivals which are recorded by a set of receivers located at a given depth. By assuming that $\Omega$ is surrounded by two homogeneous regions $\Omega_{sup}$ and $\Omega_{inf}$ (see Fig.\ref{figomega}), the receivers do not record any wave propagating below or above $\Omega$. Hence only $\Omega$ is under study. 
\begin{figure}
\begin{center}
\begin{pspicture}(10,8)
\psline{-}(0,2)(10,2)
\psline{-}(0,6)(10,6)
\uput{0.1}[90](0,6){$z=0$}
\uput{0.1}[90](0,2){$z=z_{max}$}
\uput{0.1}[0](5,7){$\Omega_{sup}$}
\uput{0.1}[0](5,4){$\Omega$}
\uput{0.1}[0](5,1){$\Omega_{inf}$}
\uput{.1}[0](8,7){$v=c^{sup}$}
\uput{.1}[0](8,4){$v=c(\vec x)$}
\uput{.1}[0](8,1){$v=c_{inf}$}

\end{pspicture}
\end{center}
\mcaption{\label{figomega} The propagation medium}
\end{figure}

The phenomenon is governed by the wave
equation set in $\rr^3$: 
\begin{equation}\label{propa1}
\left\{
\begin{array}{l}
\frac{1}{v^{2}}\partial^{2}_{t} p-div(\vec{\nabla} p)=S\mbox{ dans }\rr^3\times]0,T[ ,\\
p(0,\vec{x})=\partial_{t}p(0,\vec{x})=0\mbox{ dans }\rr^3,\\
\end{array}
\right.
\end{equation}
whose solution is the acoustic pressure $p=p(t,\vec{x})$.
The nonnegative variable $t$ denotes the time while $\vec x=^t(x,y,z)$ stands for the position vector.
In the following, $\vec x'=^t(x,y)$ designates the transverse variable. The propagation velocity $v:=v(\vec x)$ is constant-valued, equal to $c^{sup}$ in $\Omega_{sup}$ and $c_{inf}$ in $\Omega_{inf}$. In $\Omega$, $v:=c(\vec x)$ varies in all the directions. The source $S$ is a Ricker function modelling an
impulsion,
$$S(t,\vec x)=\delta(\vec x-\vec x_s)\frac{\partial^2}{\partial t^2}\left(e^{-\alpha^2(t-t_*)^2}\right)$$
where $\delta_{x_s}$ denotes the Dirac distribution at
$\vec{x}_s=^t(x_s,y_s,0)$, $ \alpha=\pi\nu_*$ and $t_*=1/\nu_*$. The constant $\nu_*$ is a given frequency.
We use standard notations for the Differential Calculus, such as:
$\partial_t$ denotes the partial derivative in the time variable,
if $\partial_i$ stands for the partial derivative with respect to
$i=x,y,z$, $div
\vec{v}=\partial_{x}v_{x}+\partial_{y}v_{y}+\partial_{z}v_z$
represents the divergence of a vector field $\vec{v}$ with
components $v_x, v_y, v_z$ and the operator $\vec\nabla$ is the
gradient defined as $\vec\nabla\varphi=
^{t}(\partial_{x}\varphi,\partial_{y}\varphi,\partial_{z}\varphi)$.
Model (\ref{propa1}) can be derived from two laws of the mechanic
of continuous media (see for instance \cite{Germain}) which
involve the velocity in the fluid $\vec{u}=
\vec{u}(t,\vec{x})=^t(u_x,u_y,u_z)$. The first concerns the
conservation of the displacement quantity and reads as: for
$i=x,y,z$,
\begin{equation}\label{conserv1}
\rho\frac{\partial ( u_{i})}{\partial
t}=div\vec{\sigma}_{i}+f_{i}
\end{equation}
where $\vec{f}=^{t}(f_{x},f_{y},f_{z})$ designates the density of
the exterior forces and $\rho$ represents the density of the medium supposed to be constant herein. The stress tensor $\vec\sigma$ describes the
behavior laws of the fluid whose lines are denoted by
$\vec\sigma_{i}$ with $i=x,y,z$. We use the notation:
\begin{equation}
div\vec\sigma_{i}=\frac{\partial\sigma_{ix}}{\partial
x}+\frac{\partial \sigma_{iy}}{\partial y} +\frac{\partial
\sigma_{iz}}{\partial z}.
\end{equation}
In the particular case of a perfect fluid, the tensor $\vec{
\sigma}$ can be written as :
$$\vec{\sigma}=-pI_{3}$$
where $I_{3}$ denotes the identity matrix of order 3. The law
(\ref{conserv1}) modifies then into :
\begin{equation}\label{conserv2}
\rho\frac{\partial \vec{u}}{\partial t}=-\vec{\nabla}
p+\vec{f}.
\end{equation}
In the framework of seismic prospection, the exterior forces are
often negligible. The source is an impulsion which is taken into
account in the second law describing the conservation of the mass:
\begin{equation}\label{conserv3}
\frac{1}{v^{2}\rho}\frac{\partial p}{\partial t}+div \vec{u}=q,
\end{equation}
where $q=q(t,\vec{x})$ represents an injection rate per mass
unity. By applying the divergence to Eq.(\ref{conserv2}) and by
deriving Eq.(\ref{conserv3}) with respect to the time, one can
eliminate the velocity $\vec{u}$ and the acoustic pressure $p$ is
solution to the second-order wave equation with right-hand-side (rhs) given by 
$S=\partial_{t}q$. Thus, we have:
$$q(t,\vec x)=\delta(\vec x-\vec x_s)\frac{\partial}{\partial t}\left(e^{-\alpha^2(t-t_*)^2}\right)$$
The wave equation (\ref{propa1}) can then be rewritten as a
first-order hyperbolic system of the form:
\begin{equation}
\left\{
\begin{array}{c}
\rho\partial _{t}\left(  \vec{u}\right)+\vec{\nabla} p =\vec{0}\mbox{ dans }]0,T[\times\rr^3, \\
\frac{1}{v^{2}\rho}\partial_{t} p +div\vec{u}=q
\mbox{ dans }]0,T[\times\rr^3,\\
p(0,\vec{x}) = 0 \mbox{ et } \vec{u}(0,\vec{x})=0\mbox{
dans }\rr^3.
\end{array}
\right.  \label{propa2}
\end{equation}
For the numerical simulations, the direction of the depth is selected as the one 
governing the sense of propagation. This approach is classical in
case of homogeneous media where the physical parameters do not
vary. By using the formalism of pseudo-differential operators, it
can be generalized to inhomogeneous media. In this work, we are
interested in system (\ref{propa2}) to which we associate a reduced
system which is the subject of the following section.
\section{Reduced system}\label{systred}

The initial system is given by (\ref{propa2}). We recall that the
initial data are supposed to be null. We propose to describe the
propagation of waves along the direction $Oz$.
 This idea was formerly exploited in \cite{bremmer} next in \cite{Corones,aki} where the problem was to compute
 solutions propagating into plane and homogeneous layers. In such a way, the problem turned into the solution of a
 system with constant coefficients and the solution was written as a
linear combination of eigenvectors of the problem.
  This
  approach was successful for the study of the propagation of plane waves into 1D homogeneous layers
 \cite{bremmer}  and the theory developed next in \cite{Corones} provides an extension to the 2D case.
 Nevertheless, in case of more complex media, the constitutive parameters $\rho$ and $c$ of the domain $\Omega$ vary
 in all the directions. This is why we need the formalism of pseudo-differential operators to generalize this method of
 modelling, just as it was formerly suggested by M.V. de Hoop in \cite{ref1}. \\

 In the following, we limit our attention to the solution in $\Omega$. We eliminate the time variable by applying a Laplace
 transform to the system and the related variable to $t$ is
 denoted by $\omega$.
 System (\ref{propa2}) then transforms into a stationary system, which reads in
 $\Omega$ as:
$$
\left\{
\begin{array}{l}
\vec {\nabla}_{\bot}\widehat{p}+i\omega \left( \rho
\widehat{{\vec u }}\right)
^{\prime }=\vec 0, \\
\partial_{z}\widehat{p}+i\omega \left( \rho \widehat{u_{z}}\right) =%
0, \\
\frac{i\omega}{\rho c^2}\widehat{p}+div_{\bot}\,\left(\widehat{\vec u}\right) ^{\prime }+\partial
_{z}\widehat{u}_{z}=\widehat{q},
\end{array}
\right.
$$
where
$\vec{\nabla}_{\bot}=^{t}\left(\partial_{x},\partial_{y}\right) $,
$\widehat{{\vec u}}^{\prime
}=^{t}\left( \widehat{u}_{x},\widehat{u}_{y}\right) $ and $div_{\bot}\widehat{{\vec u}}^{\prime }=\partial _{x}\widehat{u_{x}}%
+\partial _{y}\widehat{u_{y}}$. In the above equations,
$\widehat{\varphi }$ represents the Laplace transform of
$\varphi$.
 The time variable being suppressed, the transverse unknowns $\widehat{{\vec u}}^{\prime }$ are
 eliminated by plugging the first equation into the third. Then we get a system with a pair of
 unknowns of the form:
 \begin{equation}
\left( {\vec D}_{z}+L\right) {\vec U}={\vec F}  \label{propa3}
\end{equation}
with ${\vec U}=^{t}\left( \widehat{p},\widehat{v_{z}}\right)$, ${\vec D}%
_{z}={\mathbb I}_{2}\partial _{z}$ and ${\mathbb I}_{2}$
represents the $2\times 2$ identity. The operator $L$
is defined by: \\
$$L=\begin{pmatrix}
  0 & i\omega\rho \\
  \frac{i\omega}{c^2\rho}+div_{\bot}\left(\frac{1}{i\omega\rho}{\vec{\nabla}}_{\bot}\right) & 0
\end{pmatrix}$$
and the source $\vec F$ is given by ${\vec F}=^{t}\left(0,\widehat{q}\right).$\\
By factoring by $i\omega$, $L$ reads as $L=i\omega L^{\sharp }$
where
$$L^{\sharp }=
L^{\sharp }\left( {\vec x}^{^{\prime }},\dfrac{1}{\omega }%
\partial _{{\vec x}^{\prime }}\right)$$ is defined by:
$$L^{\sharp }=\left(
\begin{array}{ll}
0 & \rho \\
\frac{1}{c^2\rho} +\dfrac{1}{\omega ^{2}}div_{\bot}\left(\dfrac{1}{\rho }{\vec \nabla }%
_{\bot}\right) & 0
\end{array}
\right) .$$ According to the theory developed by H\"ormander
\cite{H}, Operator $L^{\sharp }$ is a pseudo-differential operator
of OPS$^0$ depending on the parameter $\omega$, and if
$\mathcal{L}^{\sharp }$ denotes its symbol, we have the following
representation: for any test-function $\vec{\varphi}$,
$$L^{\sharp }
{\vec \varphi }=\frac{1}{\left( 2\pi
\right) ^{2}}\int \int {\mathcal L}^{\sharp }\left( {\vec x}^{\prime },\frac{%
{\vec k}^{\prime }}{\omega }\right){\vec \varphi }\left( {\vec s}^{\prime
}\right) e^{-i\left( {\vec x}^{\prime}-{\vec s}^{\prime}\right) .k^{\prime
}}d{\vec s}^{\prime }d{\vec k}^{\prime } $$
where ${\vec s}^{\prime }=^{t}\left( s_{x},s_{y}\right) \in {\rr}^{2}$ and ${\vec k}^{\prime }=^{t}\left( k_{x},k_{y}\right) $ is the dual variable of $\vec x'$ such as the symbol of ${\vec \nabla }_{\bot}$ is given by $i{\vec k%
}^{^{\prime }}$.
\\
Symbol of $L^{\sharp }$ is defined by:
$${\mathcal L}^{\sharp }=\begin{pmatrix}
  0 & \rho \\
  \frac{1}{\rho c^2} -\dfrac{\left| {\vec k}^{\prime }\right| ^{2}}{\omega
^{2}\rho } & 0
\end{pmatrix}
 $$
with $\left| {\vec k}^{\prime }\right| ^{2}=k_{x}^{2}+k_{y}^{2}$.
Each term of symbol ${\cal L}^{\sharp }$ is an homogeneous
function of ${\vec k}^{\prime }$ with degree 0. Hence it belongs to
a class of matrices whose terms are in $S^0$. Then it admits an
asymptotic development with respect to the parameter
$\frac{\left|{\vec k}^{\prime}\right|}{\omega}$ which is small at
high frequency, that is when $\omega$ is large. The operator
$\mathcal{L}^{\sharp }$ can then be seen as a semi-classical
operator. Moreover the symbol of $L=i\omega L^\sharp$ admits an
asymptotic development according to the small parameter
$\omega^{-1}$. Hence we are in the framework of the
$\omega$-calculus developed by O. Lafitte in \cite{lafitte}
and $L$ can be considered as a $\omega$-pseudo-differential. This
property is interesting for the construction of high-order models
and is the subject of a current work \cite{bala}.\\
To solve (\ref{propa3}) turns into describing the propagation of
the field ${\vec U}$ along the depth $z$. The formalism of
pseudo-differential operators is well-adapted to extend the
classical way for solving differential systems with constant
coefficients and based on the diagonalization of matrix $L$. In the 
case of pseudo-differential operators (or differential operators
with variable coefficients ), Taylor \cite{T} has developed a
method of factoring strictly hyperbolic systems. This method
allows to replace the initial system by a system of two one-way
equations. One will describe the down-going propagation and the
other will model the up-going displacement. Both of the equations
are coupled by a off-diagonal matrix which are assimilated to
reflection terms. The idea of Taylor has been developed in
\cite{ref1,ref2}, later in \cite{ref3} for acoustic waves, in the
simplest case where $\rho$ is constant. Notice also that it has
been applied in \cite{AntoineBarucq} for the construction of
radiation conditions for electromagnetism.

\section{First-order Formulation} \label{formsimple}
The first-order formulation reads as a transport equation whose
derivation is based on the factoring of the symbol $ {\mathcal
L}^{\sharp }$.
 The eigenvalues of  $ {\mathcal L}^{\sharp }$ are given by the symbols $\gamma _{0}$ and $-\gamma _{0}$
 with
 \begin{equation}\label{valpro}
 \gamma _{0}=\left( \dfrac{1}{c^{2}\left( {\vec x}\right) }-\dfrac{\left| {\vec k}^{\prime }\right| ^{2}}
 {\omega ^{2}}\right) ^{1/2}.
 \end{equation}
 In order to define the square-root involved in (\ref{valpro}), the plane $\rr ^{2}$ is divided into three regions.
 The first is defined by

 $$\mathcal{H}=\left\{\vec{k}^{\prime}\in \rr^{2},
 \left|\vec k^{\prime }\right| ^{2}< \dfrac{\omega ^{2}}{c^{2}\left( {\vec
 x}\right)}\right\}.$$
 As soon as $\vec{k}^{\prime}\in \mathcal{H}$, matrix $ {\mathcal L}_{0}^{\sharp }$
 admits two eigenvalues which are single and real. This region corresponds to the one in which System (\ref{propa2}) is
  strictly hyperbolic and then $\mathcal{H}$ is called the
  hyperbolic region. It is exactly the region in which propagation occurs and the sign of each eigenvalue indicates the
  sense of propagation. The eigenvalue $\gamma _{0}$ is associated to the downgoing part of the wave i.e. the part propagating
  in the direction of increasing $z$. In the same way the eigenvalue $-\gamma _{0}$ is related to the upgoing part of the wave.
     The second region is defined by \\
  $${\mathcal E}=\left\{\vec{k}^{\prime}\in\rr^{2},
  \left| {\vec k}^{\prime }\right| ^{2}> \dfrac{\omega ^{2}}{c^{2}\left( {\vec x}\right)
  }\right\}.$$\\
  If $\vec{k}^{\prime}\in \mathcal{E}$, the eigenvalues of $ {\mathcal L}_{0}^{\sharp }$
  are purely imaginary. System (\ref{propa2}) is then elliptic in this region which is called the elliptic region.
  It defines a subset of frequencies which are not linked to the propagation. At last the third region
  is defined by ${\mathcal G}=\left\{ \left| {\vec k}^{\prime }\right| ^{2}=\dfrac{\omega ^{2}}{%
c^{2}\left( {\vec x}\right) }\right\} $ in which (\ref{propa2})
degenerates. Indeed the eigenvalue 0 is double and the problem is
neither hyperbolic nor elliptic. Region $\mathcal{G}$ is called
the glancing region and is related to grazing rays.\\

Then we have the following result:

\begin{proposition} \label{propos1}Let $m \in \mathbb{Z}$ and $ \vec{U}$ be a solution to the reduced system.
There exists at least an operator $P_{m} \in OPS^{m}$ with inverse
$Q_{-m}\in\mbox{OPS}^{-m}$ such that, if $\vec{V}=Q_{-m}\vec{U},
\vec{V}$ is solution to:
\begin{equation}
\left( \vec{D}_{z}+i\omega \Lambda _{0}\right) \vec{ V}=
\vec R_{0}\vec{V}+Q_{-m}\vec{F}\label{init}
\end{equation}
where $\Lambda _{0}\in OPS^{0}$ is a unique diagonal operator and
$\vec R_{0}\in OPS^{0}$ depends on $P_{m}$.
\end{proposition}

\begin{proof}
Let ${\mathcal M}_{0}$ be the matrix defined as
$${\mathcal M}_{0}=\begin{pmatrix}
  \gamma _{0} & 0 \\
  0 & -\gamma _{0}
\end{pmatrix}$$
We aim at constructing a matrix $\mathcal{P}$ such that $\mathcal{M}_0=\mathcal{P}^{-1}\mathcal{L}^\sharp\mathcal{P}$.
If $w=^{t}(w^{+},w^{-})$ denotes one of the
eigenvectors, associated with the eigenvalue $\pm\gamma_0$, its coordinates satisfy the equation:
$$ \rho w^{-}=\pm \gamma_{0} w^{+}$$ which admits an infinite number of solutions. We propose to solve this equation in such
a way that the coordinates of the eigenvectors are in the same
class of symbols. Hence if $S^{m},m\in \mathbb{Z}$ denotes this
class \cite{T}, we define the conjugated matrix $\mathcal{P}_m=\left[w^+, w^-\right] $ where $w^+$ and $w^-$ have been chosen in $S^{m}$. Then $\mathcal{P}_m$ defines an operator of OPS$^{m}$ denoted by $P_{m}$. Matrix $\mathcal{P}_{m}$ is invertible and
  its inverse $\mathcal{P}_{m}^{-1}$ is the principal symbol of the inverse of $P_{m}$ denoted by $Q_{-m}\in\mbox{OPS}^{-m}$. As far as the symbols are concerned,
  we have the relation
  $$\mathcal{P}_{m}^{-1}{\mathcal L}^{\sharp }\mathcal{P}_{m}={\mathcal M}_{0}.$$
  Then we define the operator $\Lambda_0$ \emph{via} its symbol by the relation:
  $$\sigma(\Lambda _{0})={\mathcal M}_{0}=\sigma_{P}(Q_{-m}L^{\sharp}P_{m}).$$
  According to the pseudo-differential theory \cite{T}, we deduce that there exists a regularizing operator $\vec R_{-1}\in$ OPS$^{-1}$
  such that
  \begin{equation}\label{eq1}
  Q_{-m}L^{\sharp}P_{m}=\Lambda_{0}+\vec R_{-1}.
  \end{equation}

  Operator $\vec R_{-1}$ is uniquely defined by its symbol whose expression is known thanks to the composition rule of
  pseudo-differential operators \cite{T}.\\
  Next let us consider $\vec{U}$ as a solution to (\ref{propa3}). We set $\vec{V}=Q_{-m}\vec{U}$.
  Then
  $\vec{V}$ satisfies:
  $$ (\vec{D}_{z}+L)P_{m}\vec{V}=\vec{F}$$
  and by composing on the left by $Q_{-m}$, we get :
 $$ Q_{-m}(\vec{D}_{z}+L)P_{m}\vec{V}=Q_{-m}\vec{F}.$$
 Moreover since
 $$\vec{D}_{z}P_{m}=P_{m}\vec{D}_{z}+\partial_{z}P_{m}$$\\
 where $\partial_{z}P_{m}$ is the pseudo-differential operator
 with symbol $\partial_{z}\mathcal{P}_{m}$, the system modifies
 into
 $$(\vec{D}_{z}+Q_{-m}LP_{m})\vec{V} +
 Q_{-m}\partial_{z}P_{m}\vec{V} = Q_{-m}\vec{F}.$$\\
 But by definition of $L$, we have :
 $$Q_{-m}LP_{m}=i\omega\left(Q_{-m}L^{\sharp}P_{m}\right)$$\\
 which implies, by taking (\ref{eq1}) into account,
 $$(\vec{D}_{z}+i\omega\Lambda _{0}+i\omega \vec R_{-1}+Q_{-m}\partial_{z}P_{m})\vec{V}
 =Q_{-m}\vec{F}.$$\\
 Proposition 2.1 is then proved by setting :
 $$
 \vec R_{0}=-i\omega\left(\vec R_{-1}\right)-Q_{-m}\partial_{z}P_{m}.$$
 \end{proof}

 \begin{remarqueit}
 Let $m$ be fixed. Then, Matrix $\mathcal{P}_{m}$
 is not single. That means that even if $m$ is fixed, we can
 construct an infinite number of models which differ all by the
 coupling operator $\vec R_{0}$. Nevertheless, it can be chosen in such a way the computational cost is a bit cut down, as seen further.
 \end{remarqueit}

 Hence the model is defined from a transport equation which involves two one-way equations of order +1. In the rhs,
  $\vec R_0$ takes the coupling terms into account and they are modelled by the off-diagonal terms of the rhs.\\
 Our aim is to account for the lateral variations of the velocity. Hence the density mass can be constant or variable in all
 the directions. In this work, we limit our attention to the case $\rho$ is
 constant.\\

The coupling operator can be split into the sum of two operators,
the first being defined by the diagonal part of $\vec R_{0}$ and
the second by the off-diagonal part. Introduce
$$\vec R_{0}^{d}=diag(\vec R_{0})\mbox{ and }
\vec R_{0}^{ad}=\vec R_{0}-\vec R_{0}^{d}.$$ By
considering this decomposition, one can construct a new model in
which the transport equation is diagonal and involves the sum of
two operators in OPS$^{1}$ and OPS$^{0}$ respectively. Then we
get:
\begin{corollaire} \label{coro1}
Let $m \in \mathbb{Z}$. There exists at least an operator $P_{m}\in
OPS^{m}$ with inverse $Q_{-m}\in OPS^{-m}$ such that if $\vec{U}$
denotes a solution to the reduced system,
$\vec{V}=Q_{-m}\vec{U}$ is solution to:\\
\begin{equation}
\left( \vec{D}_{z}+i\omega \Lambda _{0}-\vec R_{0}^{d}\right)
\vec{ V}= \vec R_{0}^{ad}\vec{V}+Q_{-m}\vec{F}
\end{equation}
\end{corollaire}

The model in \ref{coro1} consists of two one-way equations which
are coupled by the operator $\vec R_{0}^{ad}$. According to
its definition, $\vec R_{0}$ is not easy to define exactly because
it is given by an asymptotic expansion expressing its symbol. It
is obvious that we cannot use the exact symbol of $\vec R_{0}$ and
that we must be content with approximating the system by a truncated one. The truncation order must then be
defined. Since the initial reduced system is defined from an
operator $L$ in OPS$^{1}$ with a rhs in OPS$^0$, we propose to
keep this property in the one-way model. Hence we address the question
of putting the principal part $\vec R_{0,P}$ in place of $\vec R_{0}$.
Operator $\vec R_{0,P}$ is defined as the one whose symbol is $\sigma_P(\vec R_{0})$. Then we get the following result:
\begin{proposition} \label{propos3} Under the same assumptions than in Proposition
\ref{propos1}, an approximate one-way model is given by: if
$\vec{V}=Q_{-m}\vec{U}$, $\vec{V}$ is solution to
\begin{equation}
\left( \vec{D}_{z}+i\omega \Lambda _{0}-\vec R_{0,P}^{d}\right)
\vec{ V}=
\vec R_{0,P}^{ad}\vec{V}+Q_{-m}\vec{F},\label{tgsp}
\end{equation}
where the symbol of $\vec R_{0,P}$ is given by:
\begin{eqnarray*}
\sigma(\vec R_{0,P})&=&-\omega\left\{\mathcal{P}^{-1}_{m}
\left(\vec\nabla_{\vec{k}^{'}}\mathcal{L}^{\sharp}\right)\left(\vec\nabla_{\vec{x}^{'}}\mathcal{P}_{m}\right)+
\left(\vec\nabla_{\vec{k}^{'}}
\mathcal{P}^{-1}_{m}\right)\vec\nabla_{\vec{x}^{'}}\left(\mathcal{L}^{\sharp}\mathcal{P}_{m}\right)\right\}-\mathcal{P}_{m}^{-1}\partial_{z}\mathcal{P}_{m}
\end{eqnarray*}
\end{proposition}
\begin{proof}
We aim at computing the principal symbol of $R_{0}$. We use the definition:
$$\vec R_{0}=-i\omega\left(\vec R_{-1}-Q_{-m}\right)\partial_{z}P_{m}$$
with $\vec R_{-1}=Q_{-m}L^{\sharp}P_{m}-\Lambda_{0}$, which
implies that:
\begin{equation}\label{eq2}
\sigma(\vec R_{0,P})=-i\omega \sigma_{P}\left(\vec
R_{-1}+\right)-\sigma_{P}\left(Q_{-m}\partial_{z}P_{m}\right).
\end{equation}
We introduce the following notation. Let $\mathcal{A}$ and
$\mathcal{B}$ be two matrices whose terms are symbols. Let $\alpha
= ^{t}(\alpha_{x},\alpha_{y})$ be a multi-index in $\nn^2$. We set
$\vec\nabla_{\vec{k}^{'}}^{(\alpha)}
\mathcal{A}\vec\nabla_{\vec{x}^{'}}^{(\alpha)} \mathcal{B}$ the
product defined by: \\
$$\vec\nabla_{\vec{k}^{'}}^{(\alpha)}
\mathcal{A}\vec\nabla_{\vec{x}^{'}}^{(\alpha)}
\mathcal{B}=\sum_{|\beta|=|\alpha|}{\partial_{\vec{k}^{'}}^{\beta}\mathcal{A}\partial_{\vec{x}^{'}}^{\beta}\mathcal{B}}$$
where $\beta=^{t}(\beta_{x},\beta_{y})\in\nn^{2}$,
$\partial_{\vec{k}^{'}}^{\beta}=\partial_{k_{x}}^{\beta_{x}}\partial_{k_{y}}^{\beta_{y}}$
and
$\partial_{\vec{x}^{'}}^{\beta}=\partial_{x}^{\beta_{x}}\partial_{y}^{\beta^{y}}$.
\\
By applying the composition rule of pseudo-differential operators \cite{T},
we get:
\begin{equation}\label{eq3}
\sigma_{P}\left(Q_{-m}\partial_{z}P_{m}\right)=\mathcal{P}_{m}^{-1}\partial_{z}\mathcal{P}_{m}.
\end{equation}
Moreover we have
$$\sigma\left(Q_{-m}L^{\sharp}P_{m}\right)=\mathcal{P}_{m}^{-1}
\mathcal{L}^{\sharp}
\mathcal{P}_{m}-i\mathcal{P}^{-1}_{m}\left(\vec\nabla_{\vec{k}
^{'}}\mathcal{L}^{\sharp}\right)\left(\vec\nabla_{\vec{x}^{'}}
\mathcal{P} _{m} \right)-i\vec\nabla_{\vec{k}^{'}} \mathcal{P}
_{m}^{-1} \vec\nabla_{\vec{x} ^{'}}\left(\mathcal{L}^{\sharp}
\mathcal{P} _{m}\right)+ \mathcal{M} _{-2}$$ where $\mathcal{M}
_{-2}$ is the symbol of an operator in OPS$^{-2}$. Then, since
$$\mathcal{P}_{m}^{-1}
\mathcal{L}^{\sharp} \mathcal{P}_{m}=\sigma(\Lambda_{0}),$$ we can
deduce that
\begin{equation}\label{eq4}
\sigma_{P}(\vec R_{-1})=\sigma_{P}(\mathcal{P}_{m}^{-1}
\mathcal{L}^{\sharp}
\mathcal{P}_{m}-\Lambda_{0})=-i\mathcal{P}_{m}^{-1}\left(\vec\nabla_{\vec{k}
^{'}}\mathcal{L}^{\sharp}\right)\left(\vec\nabla_{\vec{x}^{'}}
\vec{P} _{m} \right)-i\vec\nabla_{\vec{k}^{'}}
\mathcal{P}^{-1}_{m} \vec\nabla_{\vec{x}
^{'}}\left(\mathcal{L}^{\sharp} \mathcal{P} _{m}\right)
\end{equation}

By plugging (\ref{eq4}) and (\ref{eq3}) into (\ref{eq2}), we
complete the proof of Proposition \ref{propos3}.
\end{proof}

\section{Setting of the numerical method}

The numerical method is now based on the solution to the approximate one-way model:
\begin{equation}
\left( \vec{D}_{z}+i\omega \Lambda _{0}-\vec R_{0,P}^{d}\right)
\vec{ V}= \left(\vec R_{0,P}^{d}+\vec R_{0,P}^{ad}\right)\vec{V}+\vec{S} \label{pet2}
 \end{equation}
where $\Lambda_0$ is the diagonal operator in OPS$^0$ whose symbol
is the diagonal matrix
$$\mathcal{M}_0=\left(\begin{array}{cc}\gamma_0&0\\0&-\gamma_0\end{array}\right).$$
Symbol $\gamma_0=\sqrt{\frac{1}{c^2}-\frac{|\vec
k'|^2}{\omega^2}}$ and its opposite $-\gamma_0$ are the
eigenvalues of $\mathcal{L}^\sharp$. The auxiliary
unknown $\vec V$ whose components are the down-going field $V_+$
and the up-going one $V_-$ is linked to $\vec U$ by the relation
$\vec U=P_0\vec V$ where $P_0$ is the operator in OPS$^0$ whose
symbol is the matrix $\mathcal{P}_{0}$. We choose $\mathcal{P}_{0}$ such as 
$$\vec
R_{0,P}^{d}=-\frac{1}{2}\Gamma_0^{-1}\partial_z\Gamma_0I_2\mbox{
and }\vec R_{0,P}^{ad}=
\frac{1}{2}\Gamma_0^{-1}\partial_z\Gamma_0J_2$$ where
$\Gamma_0\in\mbox{OPS}^0$ is defined by symbol $\gamma_0$. Then, $\mathcal{P}_{0}$ is given by $\mathcal{P}_{0}=\left(\begin{array}{cc}\rho&\rho\\\gamma_0&-\gamma_0\end{array}\right)$.  Since
$\vec R_{0,P}^{d}$ does not play a part in coupling the
down-going and up-going fields, it can be seen as a transmission
term. Introducing a parameter $\epsilon\in\{0,1\}$, we write 
\begin{equation}
\left( \vec{D}_{z}+i\omega \Lambda _{0}-\epsilon \vec R_{0,P}^{d}\right)
\vec{ V}= \left((1-\epsilon)\vec R_{0,P}^{d}+\vec R_{0,P}^{ad}\right)\vec{V}+\vec{S}. \label{pet22}
 \end{equation}
Then if $\epsilon=0$, $\vec R_{0,P}^{d}$ acts in the right-hand-side of the model (like in the model proposed by De Hoop), and if $\epsilon=1$, $\vec R_{0,P}^{d}$ is included in the left-hand-side (like in the one-way equations proposed by Zhang \emph{et al.}).
Only $\vec R_{0,P}^{ad}$
accounts for the coupling and it can be seen as the reflection
operator. In the following, we use the notations:
$$\vec R_{0,P}^{d}=T_0I_2\;\;,\;\;\vec R_{0,P}^{ad}=R_0J_2$$
with the letters $T$ to indicate the Transmission and $R$ corresponds to the Reflection.  \\
We introduce also the operator:
$$M^\epsilon=i\omega\Lambda_0-\epsilon T_0I_2.$$
To compute the solution to (\ref{pet2}) requires to invert $\vec
D_z+M^\epsilon$ and the inverse operator is denoted by
$G^\epsilon$, the so-called propagator. It governs the propagation
along the depth, in the two senses. Like $M^\epsilon$,
$G^\epsilon$ is diagonal and
\begin{equation}
G^\epsilon=\left(\begin{array}{cc}G^\epsilon_+&0\\0&G^\epsilon_-\end{array}\right).\label{pet3}
\end{equation}\\
Assume that $G^\epsilon$ is given. Then (\ref{pet2}) can be
transformed into:
$$\vec V=G^\epsilon\left((1-\epsilon)T_0I_2+R_0J_2\right)\vec V+G^\epsilon\vec S$$
and if we  introduce $K^\epsilon=G^\epsilon
\left((1-\epsilon)T_0I_2+R_0J_2\right)$, we have
\begin{equation}
(I_2-K^\epsilon)\vec V=G^\epsilon \vec S. \label{ap2}
\end{equation}
Hence the solution is computed once $(I_2-K^\epsilon)$ has been
inverted. By construction, $K^\epsilon\in\mbox{OPS}^{-1}$ because
$K^\epsilon$ arises from the composition of
$G^\epsilon\in\mbox{OPS}^{-1}$
and of
$\left((1-\epsilon)T_0I_2+R_0J_2\right)\in\mbox{OPS}^0$. Hence we
have, according to \cite{T}:
\begin{equation}
(I_2-K^\epsilon)^{-1}=\sum_{j\ge
  0}\left(K^\epsilon\right)^j. \label{pet4}
\end{equation}
 This representation corresponds to a Neumann series for the inverse of $I_2-K^\epsilon$, under the assumption
 $\|K^\epsilon\|<1$. \\
Then we use (\ref{pet4}) to expand $\vec V$ in the following form:
\begin{equation}
\left\{
\begin{array}{l}
\vec V=\sum\limits_{j\ge 0}\vec V^{(j)} \\
\vec V^0=G^\epsilon \vec S \;\mbox{ et }\;\vec V^{(j)}= K^\epsilon\vec
V^{(j-1)}\mbox{,}\,j\ge 1
\end{array}
\right .
\label{ap3}
\end{equation}

This representation of $\vec V$ is called the Bremmer series and
refers to as precursory Bremmer's works \cite{bremmer} in the simple case of a 1D model. \\
Each term $\vec V^{(j)}$ has two components denoted by $V^{(j)}_+$ and $V^{(j)}_-$ where $+$ corresponds to the
down-going part propagating in the sense of increasing $z$ and $-$ is associated to the up-going part. \\
By definition of the source $\vec F=^t(0,\widehat{q})$, the
auxiliary source $\vec S$ has \emph{a priori} two components $S_+$
et $S_-$ and the first term in the series is given by:
$$V^{(0)}_+=G^\epsilon_{+}S_+\;\mbox{ et }\;V^{(0)}_-=G^\epsilon_{-}S_-.$$
Then the second term $\vec V^{(1)}$ is obtained by the formulas:
$$\left\{\begin{array}{l}
V^{(1)}_+=K^\epsilon_{11}V^{(0)}_++K^\epsilon_{12}V^{(0)}_-\\
V^{(1)}_-=K^\epsilon_{21}V^{(0)}_++K^\epsilon_{22}V^{(0)}_-\end{array}\right.
$$
where $K^\epsilon_{lm}\;,\;1\le l,m\le2$ represent the terms of $K^\epsilon$. \\
However the numerical solution is computed into a region limited
by $z=0$ and $z=z_{max}$ and which is surrounded by two
homogeneous regions. By definition $V_-^{(0)}$ represents the propagation of
the component $S_-$ of the source which is located at the surface
$z=0$. Hence the support of $V_-^{(0)}$ is embedded in $
\Omega_{sup}$. In the numerical tests we present, we do not evaluate $V_-^{(0)}$. The algorithm is evaluated with $V_-^{(0)}=0$ which does not produce any error in the numerical results because the receivers are located at the same depth than the source. Nethertheless the omputaion of $V_-^{(0)}$ is possible and the case of a source located into $\Omega$ can be considered also.
If $K^\epsilon$ is represented from its
principal symbol:
$$\sigma_P(K^\epsilon)=\sigma_P(G^\epsilon)\left(\sigma(T_0)(1-\epsilon)I_2+\sigma(R_0)J_2\right),$$
we have $\sigma_P(K^\epsilon)=0$ in $\Omega_{sup}$.
Then since $K^\epsilon$ is a pseudo-differential operator depending continuously on the parameter
$z$, $ K^\epsilon_{l2}V_-^{(0)}=0,\;1\le l\le 2$. We then deduce that in that case, it is not useful to compute
 $V_-^{(0)}$. 
Then $\vec V^{(1)}$ reads in the simplified form as:
$$V^{(1)}_+=K^\epsilon_{11}V^{(0)}_+\;\mbox{ and }\;V^{(1)}_-=K^\epsilon_{21}V^{(0)}_+$$
while the next terms are linear combinations of the down-going and
up-going fields:
\begin{equation*}
\left\{
\begin{array}{c}
V^{(j)}_+=K^\epsilon_{11}V^{(j-1)}_++K^\epsilon_{12}V^{(j-1)}_-\\
V^{(j)}_-=K^\epsilon_{21}V^{(j-1)}_++K^\epsilon_{22}V^{(j-1)}_-.
\end{array}
\right.
\end{equation*}
with $j\ge 2$. If $\epsilon =1$, each term is simpler because 
$K^\epsilon_{11}=K^\epsilon_{22}=0$, which implies that
$V^{(1)}_+=0$ and next, using the chain rule, we get:
$$V^{(2j+1)}_+=0 \;\mbox{ et }\; V^{(2j)}_-=0\;,\;j\ge 0.$$
Hence each term is defined from either a down-going field or a up-going one which shows off the uncoupling of the
computations. This property makes think out models based on paraxial equations and a comparison between one-way models and paraxial systems is a current work \cite{badupra3}.\\

To summarize, the solution of (\ref{pet2}) involves the following
steps:
\begin{enumerate}
\item decomposition of the source: $\vec S=P_0^{-1}\vec F$
\item computation of the Bremmer terms from the propagator $G^\epsilon$, the reflection operator $R_0$ and
the transmission term $T_0$. \label{cb}
\item recomposition of the solution: $\vec U=P_0\vec V$.
\end{enumerate}
The center of the algorithm is given by the item \ref{cb}. In
order to illustrate its principle, we consider the case of a
velocity model consisting of three homogeneous layers.

\section{Bremmer series for a 2D stratified medium}
The numerical solution is based on the expansion of the fields as a
Bremmer series. We assume that $\Omega$ is defined as in Fig.
\ref{profvit}. The first layer is homogeneous with thickness
$z_1-\frac{\delta}{2}$ where $\delta$ is a small parameter and the
corresponding constant velocity is denoted by $c_1$. Next on a
thin layer with thickness $\delta$, the velocity continuously
varies from $c_1$ to the constant value $c_2$ which is the
propagation velocity in the layer with thickness
$(z_2-z_1)-\delta$. Then on a thin layer with thickness $\delta$,
the velocity continuously varies from $c_2$ to the constant value
$c_3$.

\begin{figure}
\begin{center}
\begin{pspicture}(8,9)
\psline{->}(0.8,8)(7.5,8)
\uput{.1}[0](7.5,8){$c$}
\psline{->}(1,8.2)(1,1)
\uput{.1}[225](1,1){$z$}
\psline{-}(2.5,8.1)(2.5,7.9)
\uput{.1}[90](2.5,8.1){$c_1$}
\psline{-}(4,8.1)(4,7.9)
\uput{.1}[90](4,8.1){$c_2$}
\psline{-}(5.5,8.1)(5.5,7.9)
\uput{.1}[90](5.5,8.1){$c_3$}
\psline[linestyle=dashed]{-}(0.9,6.25)(2.5,6.25)
\uput{.1}[180](0.9,6.25){$z_1-\delta/2$}
\psline[linestyle=dashed]{-}(0.9,5.75)(4,5.75)
\uput{.1}[180](0.9,5.75){$z_1+\delta/2$}
\psline[linestyle=dashed]{-}(0.9,4.75)(4,4.75)
\uput{.1}[180](0.9,4.75){$z_2-\delta/2$}
\psline[linestyle=dashed]{-}(0.9,4.25)(5.5,4.25)
\uput{.1}[180](0.9,4.25){$z_2+\delta/2$}
\psline[linestyle=dashed]{-}(0.9,2)(5.5,2)
\uput{.1}[180](0.9,2){$z_{max}$}
\psline{-}(2.5,8)(2.5,6.25)
\pscurve{-}(2.5,6.25)(2.6,6.1)(3.25,6)(3.9,5.9)(4,5.75)
\psline{-}(4,5.75)(4,4.75)
\pscurve{-}(4,4.75)(4.1,4.6)(4.75,4.5)(5.4,4.4)(5.5,4.25)
\psline{-}(5.5,4.25)(5.5,1.5)
\end{pspicture}
\end{center}
\mcaption{\label{profvit} Velocity profile for $c_1<c_2<c_3$}
\end{figure}
The interfaces between each medium are flat and they are linked to variations of the velocity
which are so important as $\delta$ is small.\\
As long as the wave propagates into a homogeneous layer, no
reflection phenomenon occurs. This is well-reproduced by the
symbol $R_0$ since it is equal to $0$ as soon as the
$z$-derivative of the velocity vanishes, as in the case of a
homogeneous medium. Then we have:
\begin{equation}
\sigma(R_0)=\sigma(T_0)=0\;\mbox{if}\;z\in\left\{z\le z_1-\frac{\delta}{2}\right\}\bigcup\left[z_1+\frac{\delta}{2},z_2-\frac{\delta}{2}\right]\bigcup\left\{z\ge z_2+\frac{\delta}{2}\right\} \label{rt0}
\end{equation}

\subsection{Series terms before any discretization}

We propose to write the terms of the Bremmer series with item $0$,
$1$ and $2$. The next terms can be straightforwardly deduced from
the term with item 2. Let us assume that if
$\varphi=^t(\varphi_+,\varphi_-)$ is a test-function depending on
$\vec x'$ and $z$, we have:
\begin{equation}
G_+^\epsilon\varphi_+=\int_0^zG_+^\epsilon(h,z)\varphi_+(h)dh\label{iad1}
\end{equation}
and
\begin{equation}
G_-^\epsilon\varphi_-=\int_{z_{max}}^zG_-^\epsilon(h,z)\varphi_-(h)dh \label{iad2}
\end{equation}
where any $G_\pm^\epsilon(h,z)$ is a pseudo-differential operator
depending continuously on $h$ and $z$. It has been proven in
\cite{prat} that this assumption is available. Then the first term
in the Bremmer series is obtained from the following arguments.
 First, just as was formerly noted, it is not useful to compute $V_-^{(0)}$. This is why we set:
$$\forall z\in\left[0,z_{max}\right],\;V^{(0)}_-(z)=0$$
Next according to its definition, $V^{(0)}_+$ reads as:
$$\forall\, 0\le z\le z_{max}\;,\;V^{(0)}_+(z)=\int_0^{z}G_+^{\epsilon}(h,z)S_+(h)dh$$
As far as the second term of the series is concerned, we have:
$$\forall 0\le z\le z_{max}\;,\;V^{(1)}_+(z)=\left(K^\epsilon_{11}V^{(0)}_+\right)(z)\;,\;V^{(1)}_-(z)=\left(K^\epsilon_{21}V^{(0)}_+\right)(z)$$
where:
$$K^\epsilon_{11}=(1-\epsilon)G_+^\epsilon T_0\;\mbox{and}\;K^\epsilon_{21}=G_-^\epsilon R_0.$$
Firstly consider the down-going term $V^{(1)}_+$. Then,
Property (\ref{rt0}) implies that:
$$ \left(T_0V_+^{(0)}\right)(z)=0\;\mbox{for all }z\in\left[0,z_1-\delta/2\right].$$
and we deduce from (\ref{iad1}) that:

$$V^{(1)}_+(z)=0\;\mbox{for}\; 0\le z\le z_1-\delta /2 $$
which gives rise to
$$V^{(1)}_+(z)=(1-\epsilon)\int_{z_1-\delta /2}^{z}G_+^{\epsilon}(h,z)(T_0V^{(0)}_+)(h)dh$$
if $z\in[z_1-\delta /2\le z\le z_1+\delta/2]$. Next using
(\ref{rt0}) again, we get:
$$ \left(T_0V_+^{(0)}\right)(z)=0\;\mbox{for all }z\in\left[z_1+\delta/2,z_2-\delta/2\right]$$
and combining this relation with (\ref{iad1}) yields:
$$V^{(1)}_+(z)=(1-\epsilon)\int_{z_1-\delta /2}^{z_1+\delta /2}G_+^{\epsilon}(h,z)(T_0V^{(0)}_+)(h)dh$$
for any $z$ into the homogeneous layer $\left[z_1+\delta/2,z_2-\delta/2\right]$.\\
Now, if $z$ belongs to $[z_2-\delta /2\le z\le z_2+\delta /2]$, we
have:
 \begin{eqnarray*}
V^{(1)}_+(z)&=&(1-\epsilon)\int_{z_1-\delta /2}^{z_1+\delta /2}G_+^{\epsilon}(h,z)(T_0V^{(0)}_+)(h)dh\\
&+&(1-\epsilon)\int_{z_2-\delta /2}^{z}G_+^{\epsilon}(h,z)(T_0V^{(0)}_+)(h)dh
\end{eqnarray*}
and in the homogeneous layer $\{z\ge z_2+\delta /2\}$
 \begin{eqnarray*}
V^{(1)}_+(z)&=&(1-\epsilon)\int_{z_1-\delta /2}^{z_1+\delta /2}G_+^{\epsilon}(h,z)(T_0V^{(0)}_+)(h)dh\\
&+&(1-\epsilon)\int_{z_2-\delta /2}^{z_2+\delta /2}G_+^{\epsilon}(h,z)(T_0V^{(0)}_+)(h)dh
\end{eqnarray*}

If we choose $\epsilon =0$, $V^{(1)}_+$ corrects $V^{(0)}_+$ in which only the downward propagation of the source is taken
into account. Hence $V^{(1)}_+$ plays the role of a corrector of the propagation by introducing the transmission effects.\\

Next when $\epsilon =1$, $V^{(1)}_+$ vanishes. In fact the transmission is included into the propagator $G^1_+$ and $V^{(0)}
_+$ directly accounts for the transmission effects at each interface. \\
Now consider the up-going term. According to (\ref{iad2}),
$$V^{(1)}_-(z)=\int^z_{z_{max}}G_-^{\epsilon}(h,z)(R_0V^{(0)}_+)(h)dh=-\int_z^{z_{max}}G_-^{\epsilon}(h,z)(R_0V^{(0)}_+)(h)dh.$$
Moreover according to (\ref{rt0}), we have:
 $$V^{(1)}_-(z)=0\;\mbox{for}\;z_{max}\ge z\ge z_2+\delta /2.$$
By applying the same reasoning than previously, we get:
 $$V^{(1)}_-(z)=-\int^{z_2+\delta /2}_{z}G_-^{\epsilon}(h,z)(R_0V^{(0)}_+)(h)dh\;\;\mbox{for}\;z\in[z_2-\delta /2,z_2+\delta /2],$$
$$V^{(1)}_-(z)=-\int_{z_2-\delta /2}^{z_2+\delta /2}G_-^{\epsilon}(h,z)R_0V^{(0)}_+(h)dh\;\;\mbox{for}\; z\in[z_1+\delta /2,z_2-\delta /2], $$
\begin{eqnarray*}
V^{(1)}_-(z)&=&-\int_{z_2-\delta /2}^{z_2+\delta /2}G_-^{\epsilon}(h,z)(R_0V^{(0)}_+)(h)dh\\
&&-\int^{z_1+\delta /2}_{z}G_-^{\epsilon}(h,z)(R_0V^{(0)}_+)(h)dh\;\;\mbox{for}\; z\in[z_1-\delta /2,z_1+\delta /2],
\end{eqnarray*}
and
\begin{eqnarray*}
V^{(1)}_-(z)&=&-\int_{z_2-\delta /2}^{z_2+\delta /2}G_-^{\epsilon}(h,z)(R_0V^{(0)}_+)(h)dh\\
&&-\int^{z_1+\delta /2}_{z_1-\delta /2}G_-^{\epsilon}(h,z)(R_0V^{(0)}_+)(h)dh\;\;\mbox{for}\; z\in[0,z_1-\delta /2].
\end{eqnarray*}

In the case where $\epsilon=0$, $G^0$ represents the propagation
without accounting for the discontinuity of the medium. Hence
$V^{(1)}_-$ results from the propagation of the reflected part of
$V^{(0)}_+$.\\
To construct $\vec V^{(2)}$ we follow the same scheme than for
$\vec V^{(1)}$:
\begin{equation*}
\left\{
\begin{array}{c}
V^{(2)}_+=K^\epsilon_{11}V^{(1)}_++K^\epsilon_{12}V^{(1)}_-\\
V^{(2)}_-=K^\epsilon_{21}V^{(1)}_++K^\epsilon_{22}V^{(1)}_-.
\end{array}
\right.
\end{equation*}
The component $V^{(2)}_+$ reads as follows:
$$\forall 0\le z\le z_{max}\;,\;V^{(2)}_+(z)=\int_0^zG_+^{\epsilon}(h,z)((1-\epsilon)(T_0V^{(1)}_+)(h)+(R_0V^{(1)}_-)(h))dh.$$
In the first layer, we have seen that $\sigma(T_0)=\sigma(R_0)=0$.
Thus we can deduce that
$$\forall\,0\le z\le z_1-\delta /2,\;\;V^{(2)}_+(z)=0.$$
Next developing the same ideas than for the construction of
$V^{(1)}_+$, we get the following relations:
\begin{description}
\item[for $z_1-\delta /2\le z\le z_1+\delta /2$] $$V^{(2)}_+(z)=\int_{z_1-\delta /2}^{z}G_+^{\epsilon}(h,z)\left((1-\epsilon)(T_0V^{(1)}_+)(h)+(R_0V^{(1)}_-)(h)\right)dh,$$
\item[for $z_1+\delta /2\le z\le z_2-\delta /2$] $$V^{(2)}_+(z)=\int_{z_1-\delta /2}^{z_1+\delta /2}G_+^{\epsilon}(h,z)\left((1-\epsilon)(T_0V^{(1)}_+)(h)+(R_0V^{(1)}_-)(h)\right)dh,$$
\item[for $z_2-\delta /2\le z\le z_2+\delta /2$] \begin{eqnarray*}
V^{(2)}_+(z)&=&\int_{z_1-\delta /2}^{z_1+\delta /2}G_+^{\epsilon}(h,z)\left((1-\epsilon)(T_0V^{(1)}_+)(h)+(R_0V^{(1)}_-)(h)\right)dh\\
&+&\int_{z_2-\delta /2}^{z}G_+^{\epsilon}(h,z)\left((1-\epsilon)(T_0V^{(1)}_+)(h)+(R_0V^{(1)}_-)(h)\right)dh,
\end{eqnarray*}
\item[and for $z\ge z_2+\delta /2$] \begin{eqnarray*}
V^{(2)}_+(z)&=&\int_{z_1-\delta /2}^{z_1+\delta /2}G_+^{\epsilon}(h,z)(R_0V^{(1)}_-)(h)dh\\
&+&\int_{z_2-\delta /2}^{z+\delta /2}G_+^{\epsilon}(h,z)\left((1-\epsilon)(T_0V^{(1)}_+)(h)+(R_0V^{(1)}_-)(h)\right)dh.
\end{eqnarray*}
\end{description}
Formulas describing $V^{(2)}_+$ into the different layers are more intricate than the ones given $V^{(1)}_+$ because of
$V^{(1)}_-\neq0$. Now the down-going term is a linear combination of down-going and up-going terms.\\

In the same way, $V^{(1)}_-$ satisfies the relation:
$$\forall 0\le z\le z_{max}\;,\;V^{(2)}_-(z)=-\int_z^{z_{max}}G_-^\epsilon(h,z)((1-\epsilon)(T_0V^{(1)}_-)(h)+(R_0V^{(1)}_+)(h))dh$$
which can be split into each layer as:
\begin{description}
\item[if $z_{max}\ge z\ge z_2+\delta /2 $] $$V^{(2)}_-(z)=0,$$
\item[if $z_2+\delta /2\ge z\ge z_2-\delta /2$] $$V^{(2)}_-(z)=-\int^{z_2+\delta /2}_{z}G_-^{\epsilon}(h,z)\left((1-\epsilon)(T_0V^{(1)}_-)(h)+(R_0V^{(1)}_+)(h)\right)dh,$$
\item[if $z_2-\delta /2\ge z\ge z_1+\delta /2$] $$V^{(2)}_-(z)=-\int_{z_2-\delta /2}^{z_2+\delta /2}G_-^{\epsilon}(h,z)\left((1-\epsilon)(T_0V^{(1)}_-)(h)+(R_0V^{(1)}_+)(h)\right)dh,$$
\item[if $z_1+\delta /2\ge z\ge z_1-\delta /2$] \begin{eqnarray*}
V^{(2)}_-(z)&=&-\int_{z_2-\delta /2}^{z_2+\delta /2}G_-^{\epsilon}(h,z)\left((1-\epsilon)(T_0V^{(1)}_-)(h)+(R_0V^{(1)}_+)(h)\right)dh\\
&&-\int^{z_1+\delta /2}_{z}G_-^{\epsilon}(h,z)\left((1-\epsilon)(T_0V^{(1)}_-)(h)+(R_0V^{(1)}_+)(h)\right)dh,
\end{eqnarray*}
\item[and if $z_1+\delta /2\ge z\ge 0$] \begin{eqnarray*}
V^{(2)}_-(z)&=&-\int_{z_2-\delta /2}^{z_2+\delta /2}G_-^{\epsilon}(h,z)\left((1-\epsilon)(T_0V^{(1)}_-)(h)+(R_0V^{(1)}_+)(h)\right)dh\\
&&-\int^{z_1+\delta /2}_{z_1-\delta /2}G_-^{\epsilon}(h,z)\left((1-\epsilon)(T_0V^{(1)}_-)(h)+(R_0V^{(1)}_+)(h)\right)dh.
\end{eqnarray*}
\end{description}
In the case where $\epsilon=1$,  $V_-^{(2)}$ vanishes since $V_+^{(1)}$ vanishes too.\\

The next terms $V_+^{(j)}$  and $V_-^{(j)}$ with $j\ge 3$ read exactly in the same way than in the case $j=2$
but by replacing the subscript $(2)$ by $(j)$ and the subscript $(1)$ by $(j-1)$ into the formulas. \\

In the case where $\epsilon=1$, the method requires two times less
operations than in the case where $\epsilon=0$. Hence one could
already suppose that it will be more judicious to choose
$\epsilon=1$. This is why we propose to investigate this question
in the following section devoted to numerical experiments.

\begin{figure}
\begin{center}
\begin{pspicture}(14,10)
\psline(1,1)(13,1)
\psline(1,3)(13,3)
\psline(1,6)(13,6)
\psline(1,9)(13,9)
\uput{.1}[180](1,6){$z_1$}
\uput{.1}[180](1,3){$z_2$}
\uput{.1}[180](1,1){$z_{max}$}
\psline[linewidth=0.05]{-}(1.5,9.1)(1.5,1)
\uput{.1}[90](1.5,9.1){Source}
\psline[linewidth=0.05,linestyle=dotted]{->}(1.5,1)(1.5,0.5)
\pscircle*(1.5,6){.08}
\pscircle*(1.5,3){.08}
\uput{.1}[270](1.5,0.5){$V_+^{(0)}=G_+^0 S_+$}
\psline[linecolor=blue]{-}(1.5,6)(4.5,9)
\psline[linecolor=blue,linestyle=dashed]{->}(4.5,9)(5,9.5)
\uput{.1}[45](1.5,6.5){\textcolor{blue}{$G_-^0(R_0V_+^{(0)})(z_1)$}}
\uput{.1}[0](5,9.5){\textcolor{blue}{$V_-^{(1)}$}}
\psline[linecolor=blue]{-}(1.5,3)(7.5,9)
\psline[linecolor=blue,linestyle=dashed]{->}(7.5,9)(8,9.5)
\uput{.1}[45](5.5,6.5){\textcolor{blue}{$G_+^0(T_0V_+^{(0)})(z_2)$}}
\uput{.1}[0](8,9.5){\textcolor{blue}{$V_-^{(1)}$}}
\psline[linecolor=red]{-}(1.5,6)(6.5,1)
\psline[linecolor=red,linestyle=dashed]{->}(6.5,1)(7,0.5)
\uput{.1}[45](3.5,4){\textcolor{red}{$G_+^0(T_0V_+^{(0)})(z_1)$}}
\uput{.1}[0](7,0.5){\textcolor{red}{$V_+^{(1)}$}}
\psline[linecolor=red]{-}(1.5,3)(3.5,1)
\psline[linecolor=red,linestyle=dashed]{->}(3.5,1)(4,0.5)
\uput{.1}[45](2,2){\textcolor{red}{$G_-^0(R_0V_+^{(0)})(z_2)$}}
\uput{.1}[0](4,0.5){\textcolor{red}{$V_+^{(1)}$}}
\end{pspicture}
\end{center}
\begin{center}
\begin{pspicture}(14,10)
\psline(1,1)(13,1)
\psline(1,3)(13,3)
\psline(1,6)(13,6)
\psline(1,9)(13,9)
\uput{.1}[180](1,6){$z_1$}
\uput{.1}[180](1,3){$z_2$}
\uput{.1}[180](1,1){$z_{max}$}
\psline[linewidth=0.05]{-}(1.5,9.1)(1.5,1)
\uput{.1}[90](1.5,9.1){Source}
\psline[linewidth=0.05,linestyle=dotted]{->}(1.5,1)(1.5,0.5)
\pscircle*(1.5,6){.08}
\pscircle*(1.5,3){.08}
\uput{.1}[270](1.5,0.5){$V_+^{(0)}=G_+^0 S_+$}
\psline[linecolor=blue]{-}(1.5,6)(4.5,9)
\psline[linecolor=blue,linestyle=dashed]{->}(4.5,9)(5,9.5)
\uput{.1}[45](1.5,6.5){\textcolor{blue}{$G_-^0(R_0V_+^{(0)})(z_1)$}}
\uput{.1}[0](5,9.5){\textcolor{blue}{$V_-^{(1)}$}}
\psline[linecolor=blue]{-}(1.5,3)(7.5,9)
\psline[linecolor=blue,linestyle=dashed]{->}(7.5,9)(8,9.5)
\uput{.1}[45](5.5,6.5){\textcolor{blue}{$G_+^0(T_0V_+^{(0)})(z_2)$}}
\uput{.1}[0](8,9.5){\textcolor{blue}{$V_-^{(1)}$}}
\pscircle[fillstyle=solid,fillcolor=blue](4.5,6){.08}
\psline[linecolor=red]{-}(1.5,6)(6.5,1)
\psline[linecolor=red,linestyle=dashed]{->}(6.5,1)(7,0.5)
\uput{.1}[45](3.5,4){\textcolor{red}{$G_+^0(T_0V_+^{(0)})(z_1)$}}
\uput{.1}[0](7,0.5){\textcolor{red}{$V_+^{(1)}$}}
\psline[linecolor=red]{-}(1.5,3)(3.5,1)
\psline[linecolor=red,linestyle=dashed]{->}(3.5,1)(4,0.5)
\uput{.1}[45](2,2){\textcolor{red}{$G_-^0(R_0V_+^{(0)})(z_2)$}}
\uput{.1}[0](4,0.5){\textcolor{red}{$V_+^{(1)}$}}
\pscircle[fillstyle=solid,fillcolor=red](4.5,3){.08}
\psline[linecolor=green]{-}(4.5,6)(9.5,1)
\psline[linecolor=green,linestyle=dashed]{->}(9.5,1)(10,0.5)
\uput{.1}[0](7,3.5){\textcolor{green}{$G_+^0(R_0G_-^0(R_0V_+^{(0)})(z_2))(z_1)$}}
\uput{.1}[0](10,0.5){\textcolor{green}{$V_+^{(2)}$}}
\psline[linecolor=green]{-}(4.5,3)(8,1)
\psline[linecolor=green,linestyle=dashed]{->}(8,1)(8.5,0.714285)
\uput{.1}[270](8,2){\textcolor{green}{$G_+^0(T_0G_+^0(T_0V_+^{(0)})(z_1))(z_2)$}}
\uput{.1}[0](8.5,0.5){\textcolor{green}{$V_+^{(2)}$}}
\pscircle[fillstyle=solid,fillcolor=green](7.5,3){.08}
\psline[linecolor=yellow]{-}(4.5,6)(9,9)
\psline[linecolor=yellow,linestyle=dashed]{->}(9,9)(9.5,9.333333)
\uput{.1}[45](7.5,8){\textcolor{yellow}{$G_-^0(T_0(G_-^0(R_0V_+^{(0)})(z_2)))(z_1)$}}
\uput{.1}[0](9.5,9.5){\textcolor{yellow}{$V_-^{(2)}$}}
\psline[linecolor=yellow]{-}(4.5,3)(10.5,9)
\psline[linecolor=yellow,linestyle=dashed]{->}(10.5,9)(11,9.5)
\uput{.1}[0](7,5.5){\textcolor{yellow}{$G_-^0(R_0(G_+^0(T_0V_+^{(0)})(z_1)))(z_2)$}}
\uput{.1}[0](11,9.5){\textcolor{yellow}{$V_-^{(2)}$}}
\pscircle[fillstyle=solid,fillcolor=yellow](7.5,6){.08}
\end{pspicture}
\end{center}
 \mcaption{\label{aaahh}The Bremmer terms}
\end{figure}

\subsection{Numerical Bremmer terms}
Let $\Delta z$ be the depth step. In practise, $\Delta z=\delta$ where $\delta$ is the thickness of the layer in which
the velocity varies continuously between two constant values. Each of the integrals involved into the Bremmer terms is
approximated by the rectangle method and to choose $\Delta z=\delta$ implies that the integrals in the variable $z$ is replaced by the multiplication by $\Delta z$.\\
For instance, when the receivers are located at the surface $z=0$, they are supposed to record
$V_-(0)$ which is obtained step by step by computing the terms of
the Bremmer series and we propose to truncate the series up to the
third term. The approximate field is then given by:

$$V_{-,app}(0)= V^{(0)}_-(0)+V^{(1)}_-(0)+V^{(2)}_-(0)+V^{(3)}_-(0)$$
with
\begin{eqnarray*}
V^{(0)}_-(0)&=&0,\\
V^{(1)}_-(0)&=&-\Delta z\left(G_-^\epsilon(z_2,0)R_0V^{(0)}_+(z_2)+G_-^\epsilon(z_1,0)R_0V^{(0)}_+(z_1)\right),\\
V^{(2)}_-(0)&=&-\Delta z\left(G_-^\epsilon(z_2,0)R_0V^{(1)}_+(z_2)+(1-\epsilon)G_-^\epsilon(z_1,0)T_0V^{(1)}_-(z_1)\right),\\
V^{(3)}_-(0)&=&-\Delta z\left(G_-^\epsilon(z_2,0)R_0V^{(2)}_+(z_2)+(1-\epsilon)G_-^\epsilon(z_1,0)T_0V^{(2)}_-(z_1)\right).\\
\end{eqnarray*}
The terms are obtained by following always the same scheme. We focus our discussion on $V_-^{(1)}$ computed at a given depth
$z$. Nevertheless the scheme depicted at Fig. \ref{aaahh} gives a general survey of the contribution of each term
of the series.\\
Assuming the propagator $G_\pm^\epsilon(.,.)$ satisfies the
Chasles relation:
\begin{equation}
G_\pm^\epsilon(z,h)G_\pm^\epsilon(h,z^\ast)=G_\pm^\epsilon(z,z^\ast)=G_\pm^\epsilon(h,z^\ast)G_\pm^\epsilon(z,h), \label{chasle}
\end{equation}
the approximate term can be written as a function of the source:

\begin{eqnarray*}
&&V_{-,app}(0)=\underbrace{-\Delta z\left(G_-^\epsilon(z_1,0)R_0G_+^\epsilon(0,z_1)S_+\right) }_{\mbox{Reflection at the first interface}}\\
&&\\
&-&\underbrace{\Delta z\left(G_-^\epsilon(z_1,0)(Id-\Delta z(1-\epsilon)T_0)G_-^\epsilon(z_2,z_1)R_0G_+^\epsilon(z_1,z_2)(Id+\Delta z(1-\epsilon))T_0G_+^\epsilon(0,z_1)S_+\right)}_{\mbox{Reflection at the second interface}}\\
&&\\
&+&\underbrace{(\Delta z)^3G_-^\epsilon(z_2,0)R_0G_+^\epsilon(z_1,z_2)R_0G_-^\epsilon(z_2,z_1)R_0G_+^\epsilon(0,z_2)S_+ }_{\mbox{First multiple}}
\end{eqnarray*}
In Geophysics, the two first terms are the most important. But it
is also of great outstanding to compute the multiples also because in the
case where the velocity varies strongly, they can have a large
amplitude and a simultaneous calculus of the two first events with
the multiples can give rise to polluted results which are very
difficult to interpret correctly. Hence it is very interesting to
dispose of a numerical method which allows to compute the
multiples and to uncouple
them from the first events. Then the multiples can be removed from the primary fields which facilitates the interpretation
 of the results. To compute the multiples, it is necessary to combine the reflection operator at least three times.
  Precisely the $m$-th multiple is the result of $2m+1$ reflections. \\
In order to explain the formula giving $V_{-,app}(0)$, we restrict
ourselves to write $V_-^{(1)}(z)$ at the depth $z$. The other
terms are obtained in the same way. Before discretization,
$V_-^{(1)}(z)$ is given by:
$$V_-^{(1)}(z)=\int_{z_{max}}^z G_-^\epsilon(h,z)(R_0V_+^{(0)})(h)dh.$$
According to the Chasles (\ref{chasle}) relation, we have:
$$\forall\;z^*\in [z,z_{max}],\;G_-^\epsilon(h,z)=G_-^\epsilon(z^*,z)G_-^\epsilon(h,z^*)$$
which implies
$$V_-^{(1)}(z)=G_-^\epsilon(z^*,z)\int_{z_{max}}^{z^*} G_-^\epsilon(h,z^*)(R_0V_+^{(0)})(h)dh+\int_{z^*}^z G_-^\epsilon(h,z)(R_0V_+^{(0)})(h)dh$$
Thus we get the relation: for any $z^*\in [z,z_{max}]$,
$$V_-^{(1)}(z)=G_-^\epsilon(z^*,z)V_-^{(1)}(z^*)+\int_{z^*}^z G_-^\epsilon(h,z)(R_0V_+^{(0)})(h)dh$$
In practise, $z^*=z+\Delta z$. Hence if $V_-^{(1)}(z)$ denotes the
discrete value of $V_-^{(1)}$, we have by approximating the second
integral by a rectangle method,
$$V_-^{(1)}(z)=G_-^\epsilon(z+\Delta z,z)V_-^{(1)}(z+\Delta z)-\Delta zG_-^\epsilon(z+\Delta z,z)(R_0V_+^{(0)})(z+\Delta z). $$\\

The case of a stratified medium is very useful to show clearly why the cases $\epsilon=0$ and $\epsilon=1$ are different. Indeed, since the symbols do not depend on the variable $\vec x'$, we can write the propagator in the wave number domain and the symbol $t_0$ of $T_0$ is non zero only at each interface. Hence if we consider a layer $\{z_1<z<z_2\}$, the propagators $G_\pm^\epsilon$ are represented in the Fourier variables as \cite{ref1, ref3, prat}:
$$g_+^\epsilon=e^{-i(z_2-z_1)\omega\gamma_0}e^{\epsilon\Delta zt_0(z_1)}$$
and
$$g_-^\epsilon=e^{-i(z_2-z_1)\omega\gamma_0}e^{-\epsilon\Delta zt_0(z_1)}$$
By injecting these relations into the definition of $V_{-,app}(0)$, we get for the primary reflections:
\begin{itemize}
\item when $\epsilon=0$
\end{itemize}
\begin{eqnarray*}
&\widehat{V_{-,app}^P}(0)=e^{-iz_1\omega\gamma_0}(-\Delta zr_0(z_1))e^{-iz_1\omega\gamma_0}\widehat{S}_+&\\
&+e^{-iz_1\omega\gamma_0}(Id-(\Delta zt_0(z_1)))e^{-i(z_2-z_1)\omega\gamma_0}(-\Delta zr_0(z_2))e^{-i(z_2-z_1)\omega\gamma_0}(Id+(\Delta zt_0(z_1)))e^{-iz_1\omega\gamma_0}\widehat{S}_+&\\
\end{eqnarray*}
\begin{itemize}
\item and when $\epsilon=1$
\end{itemize}
\begin{eqnarray*}
&&\widehat{V_{-,app}^P}(0)=e^{-iz_1\omega\gamma_0}(-\Delta zr_0(z_1))e^{-iz_1\omega\gamma_0}\widehat{S}_+\\
&+&e^{-iz_1\omega\gamma_0}e^{-\Delta zt_0(z_1)}e^{-i(z_2-z_1)\omega\gamma_0}(-\Delta zr_0(z_2))e^{-i(z_2-z_1)\omega\gamma_0}e^{-\Delta zt_0(z_1)}e^{-iz_1\omega\gamma_0}\widehat{S}_+\\
\end{eqnarray*}
We can then observe that the term $(Id\mp(\Delta zt_0(z_1)))$ acting in the case $\epsilon=0$ is replaced by the exponential $e^{\mp\Delta zt_0(z_1)}$ in the case $\epsilon=1$. These two terms are of the same order when $\Delta zt_0(z_1)$ is small, which is the case when the medium is smooth, \emph{i.e.} when the vertical variations of the velocity are small.

\section{Description of the software} \label{pgm}

We use the language Fortran 90. The code has been developed from the  C one written by J\'er\^ome Le Rousseau. It runs
in the 2D and 3D cases. Its architecture has been modified as compared to the first version. Each routine is now
independent of the other ones and can be replaced by anyone else.  This is an interesting property for the analysis of the numerical method which
can now be coupled with other methods like the ones involving finite difference schemes for instance.
The computation of the Bremmer series has been optimized in the case where $\epsilon=0$ \cite{prat}. Hence, the computational cost is almost the same for $\epsilon=0$ and $\epsilon=1$.

\subsection*{\underline{Data}} The source $F$ is obtained from a FFT (Fast Fourier Transform) of the source $q$.
One gets a table as a function of the frequency $\omega$. The
source $q$ is a Ricker function which implies that its FFT
decreases quickly to zero. This is why we define a window of
calculus into the interval $[0,\omega_{max}]$. Next one sets the
number of terms in the Bremmer series by pointing at the number $N$
of multiples which must be modelled.  Then the code will compute $V_-^{(j)}$ up to $j=2N+1$.
\subsection*{\underline{Decomposition of $\vec F$}}
The cases $\epsilon=0$ and $\epsilon=1$ are distinguished. In both cases, $\omega$ varies from $0$ to $\omega_{max}$ and

\begin{enumerate}[{(}i{)}]
\large
\item \underline{$\epsilon=1$}
\normalsize
We know that when  $\epsilon=1$, $V_+^{(2j+1)}=0$ and $V_-^{(2j)}=0$.\\
The integer $m$ varies from $0$ to $N$ and one computes
$V_+^{(2m)}$ et $V_-^{(2m+1)}$ at each depth $l\Delta z$ where $l$
varies from $0$ to $l_{max}$ with $z_{max}=l_{max}\Delta z$. This
can be summarized by:
\begin{itemize}
\item[]\underline{Initialization}
\begin{itemize}
\item[]for $l$ varying from $1$ to $l_{max}$
\begin{itemize}
\item[ ]compute $V_+^{(0)}(l\Delta z)=G_+^1\left(S_+((l-1)\Delta z)\right)$
\item[ ]store $R_0V_+^{(0)}$
\end{itemize}
\item[]stop the variations of $l$
\item[]for $l$ varying from $l_{max}-1$ to $0$
\begin{itemize}
\item[ ]compute $V_-^{(1)}(l\Delta z)=G_-^1\left(R_0V_+^{(0)}((l+1)\Delta z)\right)$
\item[ ]store $R_0V_-^{(1)}$
\end{itemize}
\item[]stop the variations of $l$
\end{itemize}
\item[]\underline{for $m$ varying from $1$ to $N$},
\begin{itemize}
\item[]for $l$ varying from $1$ to $l_{max}$
\begin{itemize}
\item[ ]compute $V_+^{(2m)}(l\Delta z)=G_+^1\left(R_0V_-^{(2m-1)}((l-1)\Delta z)\right)$
\item[ ]store $R_0V_+^{(2m)}$
\end{itemize}
\item[]stop the variations of $l$
\item[]for $l$ varying from $l_{max}-1$ to $0$
\begin{itemize}
\item[ ]compute $V_-^{(2m+1)}(l\Delta z)=G_-^1\left(R_0V_+^{(2m)}((l+1)\Delta z)\right)$
\item[ ]store $R_0V_-^{(2m+1)}$
\end{itemize}
\item[]stop the variations of $l$
\item[]The results are stored in a table $RESULT:=RESULT(\vec x',m,\omega)$.
\end{itemize}
\item[]stop the variations of $m$
\end{itemize}
\large
\item \underline{$\epsilon=0$}
\normalsize
\emph{A priori} it is the most complicated case because the downward and upward parts of each term are coupled.
 However the algorithm  proposed initially by J. Le Rousseau can be
 optimized. The optimization is based on a summation process allowing to compute both the $2j$-th and the $2j+1$-th terms
 in the same time. This reduces considerably the computational time when $\epsilon =0$ and makes it competitive with the
 case $\epsilon=1$.
\subsection*{\underline{Computation of the Bremmer terms}}
Let $j$ be the summation index into the series. We use the temporary table $TEMP:=TEMP(\vec x')$ and
the results are stored into the tables $TAB_+$ and $TAB_-$ with length $l_{max}$. \\
\underline{for $m$ varying from $0$ to $N$}
\begin{itemize}
\item[]Initialization of the downgoing part
\begin{itemize}
\item[ ]$TAEMP:=0$
\item[ ]$TAB_-(0):=0$
\end{itemize}
\item[ ]for $l=1$ to $l_{max}$,
\begin{itemize}
\item[ ]if $m=0$, $TAB_+(l-1):=0$
\item[ ]$TEMP:=G_+^0((l-1)\Delta z,l\Delta z)(TEMP-TAB_-(l-1)+TAB_+(l-1))$
\item[ ]$TAB_-(l):=-\Delta zR_0TAMP$
\end{itemize}
\item[]End of variations of $l$
\item[]Initialization of the upgoing part
\begin{itemize}
\item[ ]$TEMP:=0$
\item[ ]$TAB_+(l_{max}):=0$
\end{itemize}
\item[]for $l=l_{max}$ to $1$,
\begin{itemize}
\item[ ]$TAB_+(l):=\Delta zR_0TAMP$
\item[ ]$TEMP:=G_-^0(l\Delta z,(l-1)\Delta z)(TEMP-TAB_+(l)+TAB_-(l))$
\end{itemize}
\item[]End of variation of $l$
\item[]The final values of $TEMP$ give the values of the $m$-th multiple which are stored in the table
$RESULT:=RESULT(\vec x',m,\omega)$.
\end{itemize}
\underline{End of variations of $m$}
\end{enumerate}

\subsection*{\underline{Composition}}
For each value of $m$, we have computed the field recorded at
$z=0$ and this has been done for each value of the frequency from
$0$ to $\omega_{max}$. Then applying $P_0$ to this field, we get
the unknown $\vec U$. Next by applying an inverse FFT to its first
component, we get the value of the acoustic pressure as a function
of time. By depicting the acoustic pressure at the surface $z=0$
we get a seismic section which represents the time value of the
pressure as a function of the transverse variable $\vec x'$. In
all the pictures we state the variable $\vec x'$ is the abscissa
and the ordinate gives the time oriented to the bottom. In this
way we can represent the kinematic of the phenomenon. To depict
the dynamic, we use a grey scale for the amplitude of the acoustic
pressure at point $(\vec x',t)$ and $z=0$.\\
It is also possible to represent snapshots (value of the wave
field at $(\vec x',z,t)$). Then we use an auxiliary table
$TEMP:=TEMP(\vec x',l,\omega)$ in which we store the values of the
field at $(\vec x',z,\omega)$ and by applying an inverse FFT, we
get the field evaluated at $(\vec x',z,t)$.
\newsavebox{\organigramme}
\savebox{\organigramme}[15cm]{\footnotesize
\begin{tabular}{c}
\Ovalbox{FFT in time of the source $\vec F$}\\
 \vline \\
\Ovalbox{ \colorbox{white}{Loop on the frequences: operations representation by FFT}}\\
 \vline \\
\Ovalbox{Decomposition of the source $\vec S$}\\
 \vline \\
\Ovalbox{\colorbox{white}{ Loop on the multiples number $m$}}\\
 \vline \\
\Ovalbox{ \colorbox{white}{Loop on $1\le l\le l_{max}\}$}}\\
 \vline \\
\Ovalbox{ Propagation from $l-1$ to $l$ of the downgoing field $V_+^{(2m)}$ }\\
 \vline \\
\Ovalbox{ Reflection: couplage with the upgoing field $V_-^{(2m+1)}$}\\
 \vline \\
\Ovalbox{ Update of the field with the reflection part computed at the step  $(m-1)$}\\
 \vline \\
\Ovalbox{ Transmission: couplage with the downgoing field $V_+^{(2m)}$}\\
 \vline \\
\Ovalbox{ \colorbox{white}{End of the loop on $l$}}\\
 \vline \\
\Ovalbox{ End of the downgoing part: champ$=0$}\\
 \vline \\
\Ovalbox{ \colorbox{white}{Loop on $l_{max}\ge l\ge 0]$}}\\
 \vline \\
\Ovalbox{ Reflection: couplage with the downgoing field $V_+^{(2m+2)}$}\\
 \vline \\
\Ovalbox{ Transmission: couplage with the upgoing field $V_-^{(2m+1)}$}\\
 \vline \\
\Ovalbox{ Update of the field with the reflection part computed at the downgoing step}\\
 \vline \\
\Ovalbox{ Propagation from $l$ to $l-1 $ of the upgoing field $V_-^{(2m+1)}$ }\\
 \vline \\
\Ovalbox{ \colorbox{white}{End of the loop on $l$}}\\
\vline \\
\Ovalbox{ Record the data at the receivers}\\
 \vline \\
\Ovalbox{\colorbox{white}{End of loop on the multiples}}\\
 \vline \\
\Ovalbox{ \colorbox{white}{End of th loop on the frequencies}}\\
 \vline \\
\Ovalbox{ composition: field $\times \rho$}\\
 \vline \\
\Ovalbox{inverse FFT: seismic data }\normalsize
\end{tabular}
}
\setlength{\unitlength}{1cm}
\begin{figure}
\begin{pspicture}(15,23)
\psline{-}(7.5,1.2)(15,1.2)(15,21.4)(7.5,21.4)
\psline[linewidth=0.08]{->}(15,10.6)(15,10.8)
\psline{-}(7.5,2.4)(14.5,2.4)(14.5,19.1)(7.5,19.1)
\psline[linewidth=0.08]{->}(14.5,10.6)(14.5,10.8)
\psline{-}(7.5,4.5)(14,4.5)(14,10.1)(7.5,10.1)
\psline[linewidth=0.08]{->}(14,6.3)(14,6.5)
\psline{-}(7.5,12.2)(14,12.2)(14,17.8)(7.5,17.8)
\psline[linewidth=0.08]{->}(14,14.7)(14,14.9)
\rput[l](0,10.8){\usebox{\organigramme}}
\end{pspicture}
\nonumber\label{orga}
\end{figure}

\section{Setting of the velocity model}
\begin{figure}
 \begin{center}
\includegraphics[scale=0.6]{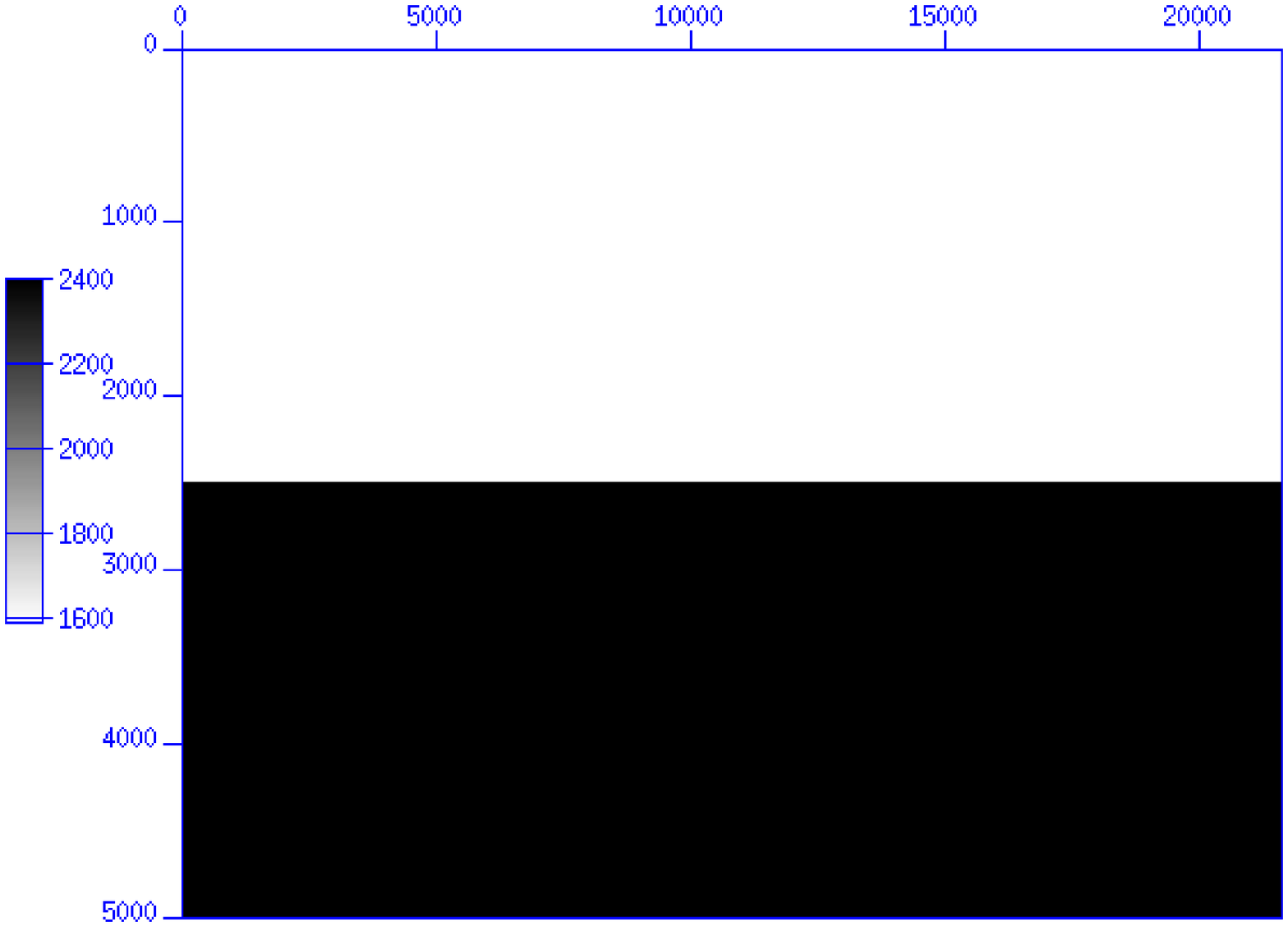}
\end{center}
\mcaption{A two layers model.}
\label{horiz}
\end{figure}

We consider a simple velocity model (See Fig.\ref{horiz}) which
consists of two layers. The topographic data of the model are the
following:
\begin{center}
\vspace{0.5cm}
\begin{tabular}{|l|c|}
\hline
depth & 5 km\\
\hline
thickness & 21.59 km\\
\hline
\end{tabular}
\vspace{0.5cm}
\end{center}
The shot is set at the middle of the domain \emph{i.e.} at $x_s=(10800,0)$ and the receivers are located at the depth $z=3000m$. By this way, we only account for the transmission effects. The transverse dimension of the model has been chosen such that
the cardinal of the discrete set of points $x_j$ is matched for
using FFTs. Each layer corresponds to a constant value of the
velocity and their interface is located at $z=2500$m. The mesh
dimensions are collected in the following table:
\begin{center}
\vspace{0.5cm}
\begin{tabular}{|l|c|}
\hline
width of the computational domain & 21.59 km \\
\hline
transverse step $\Delta x$ & 10 m \\
\hline
depth step $\Delta z$ & 10 m \\
\hline
\end{tabular}
\vspace{0.5cm}
\end{center}
In the following we refer to Velocity Model 1 (VM 1) when the velocity is equal
to 1600 ms$^{-1}$ in the first layer and 2400 ms$^{-1}$ in the
second one. VM 2 corresponds to the case where the velocity is
equal to 1600 ms$^{-1}$ in the first layer and 5000 ms$^{-1}$ in
the second one. The last model VM 3 is defined by a velocity
equal to 2400 ms$^{-1}$ in the first layer and 1600 ms$^{-1}$ in
the second one. VMs 1 and 3 are examples of velocities with a quite
low contrast while VM 2 gives an example of high velocity
contrast. \\
\begin{figure}
\begin{center}
\includegraphics[scale=0.35]{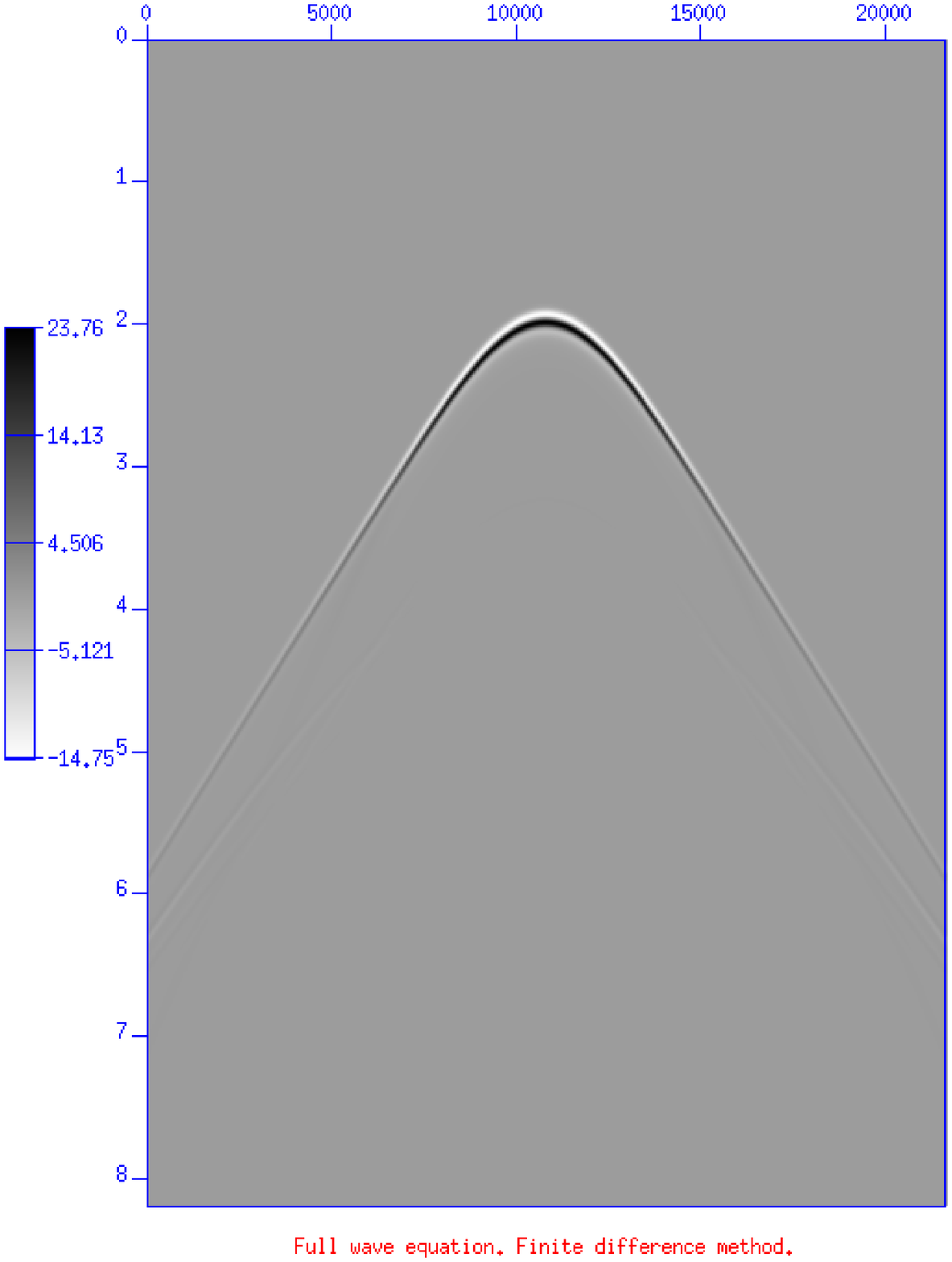}
\includegraphics[scale=0.35]{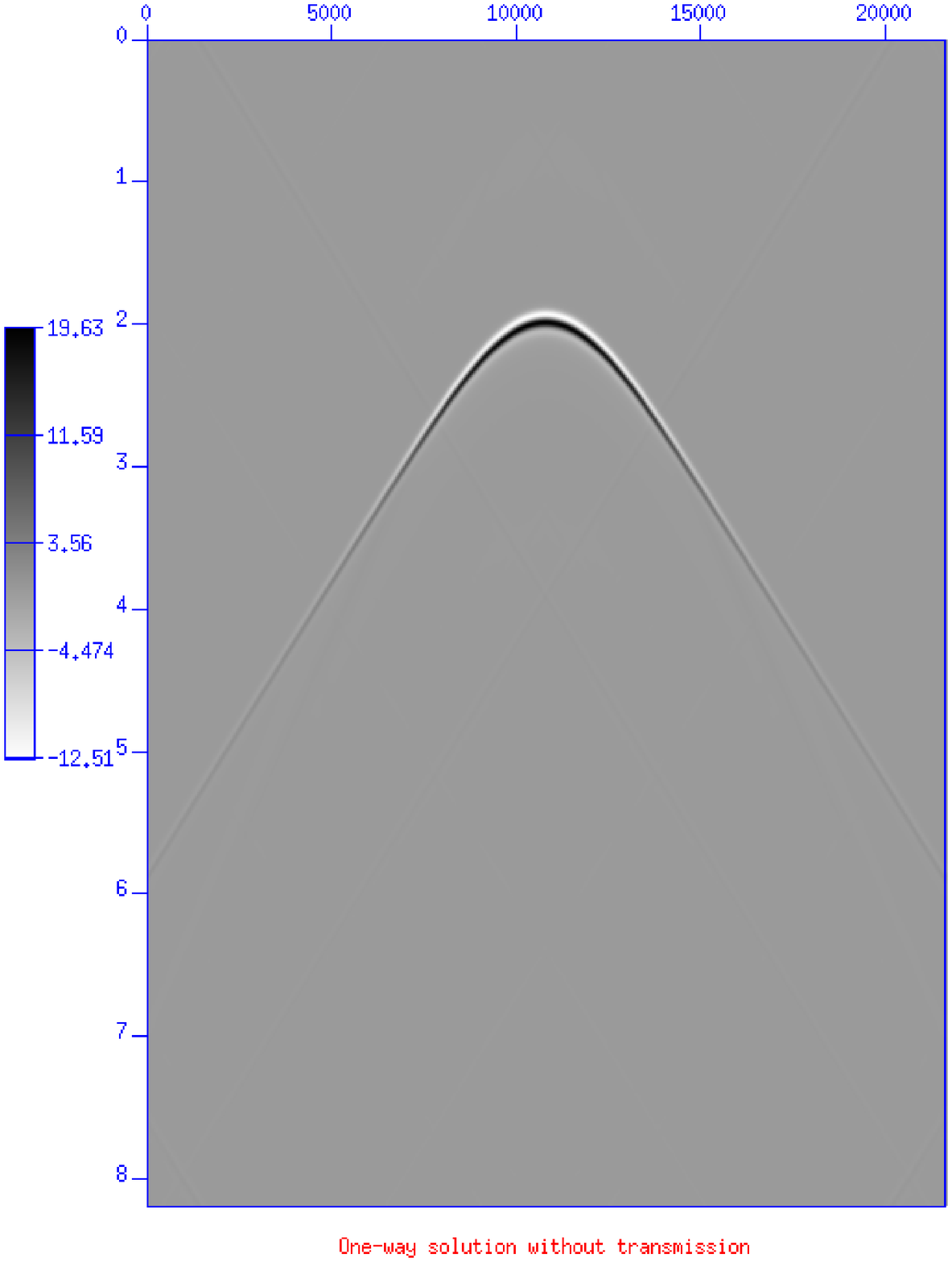}
\includegraphics[scale=0.35]{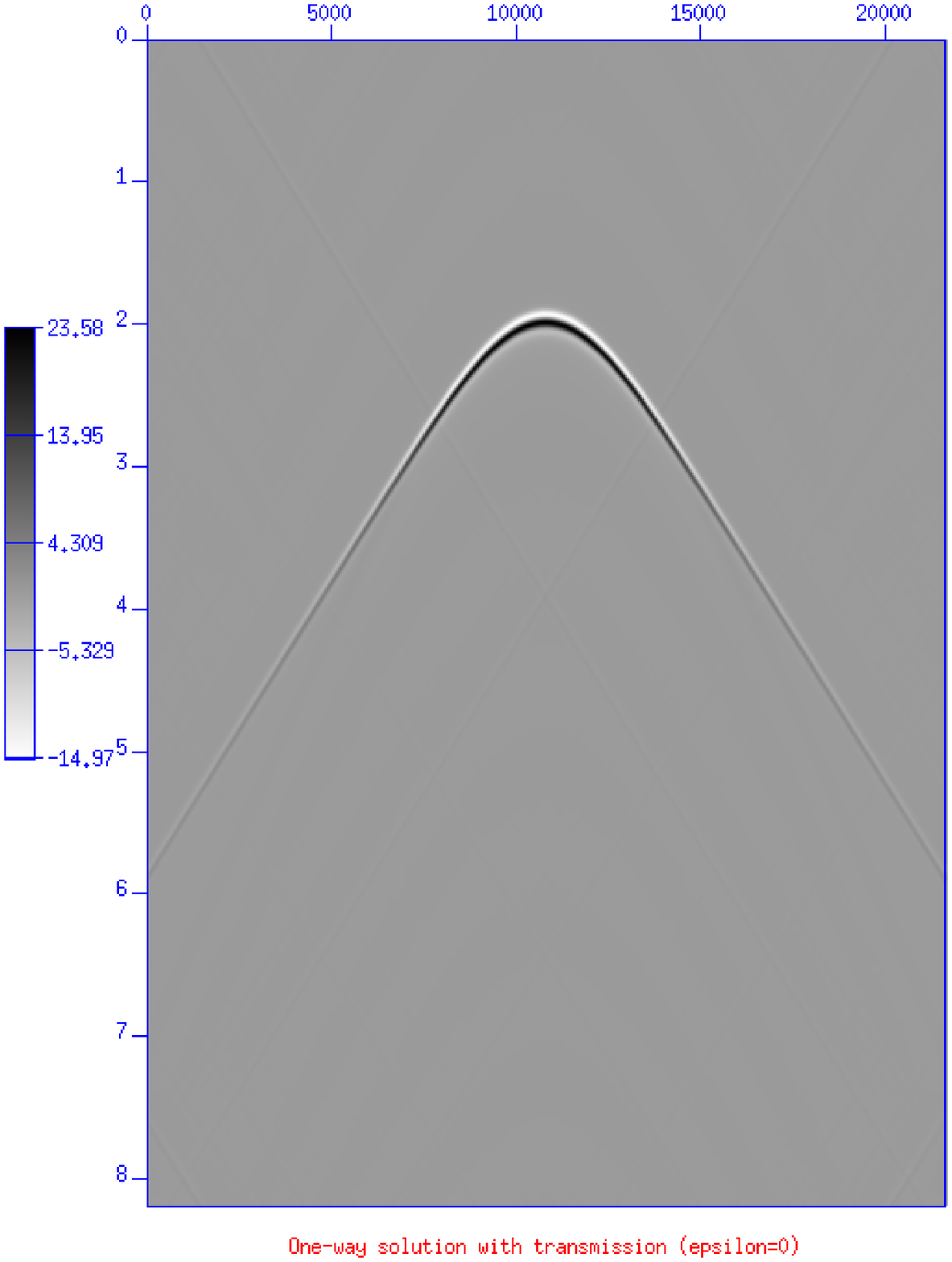}
\includegraphics[scale=0.35]{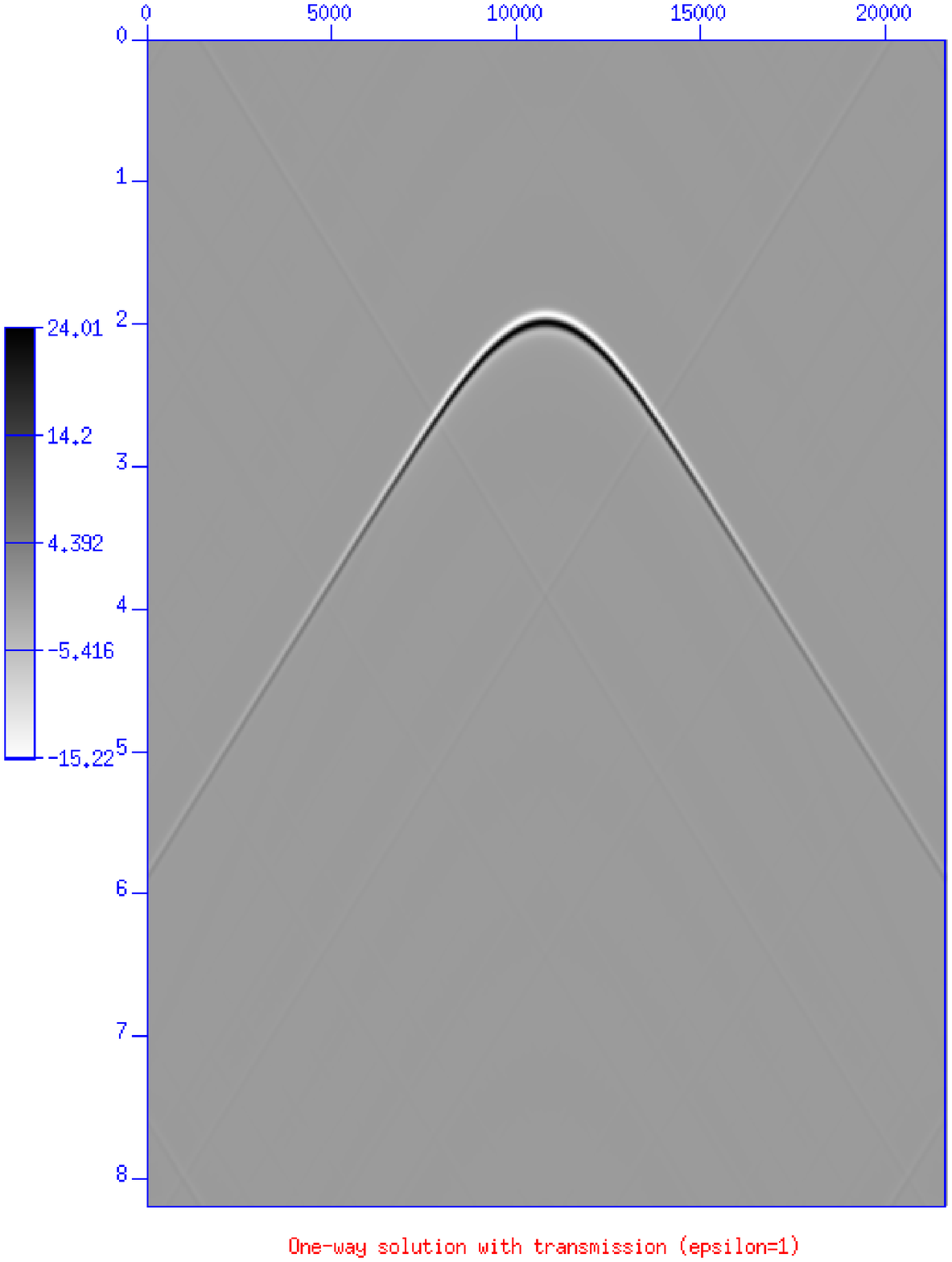}
\end{center}
\mcaption{Seismograms.}
\label{sismo}
\end{figure}
On Fig.\ref{sismo}, we have collected four seismograms. The first one depicts the reference solution obtained from the second-order full wave equation. On the right of this picture, we have set the solution of the one-way equations system as it was computed initially in \cite{ref3}. One can observe that the kinematics is well-reproduced. Nevertheless the extrema of the grey scale shows that the amplitude of the one-way solution is erroneous. This motivates the two following pictures underneath. The left seismogram depicts the best results: both the kinematics and the amplitude of the wave field are correct. The right seismogram is not so bad even if the amplitude is more erroneous than in the previous case. Hence this first collection of numerical tests shows that to include the transmission term into the model really improves the amplitude of the wave field. Moreover if we limit our analysis to consider seismograms, we can conclude that the transmission term can be numerically handled as a term of the right-hand-side of the one-way system or as a proper term of the one-way equations. But the seismograms give a global view of the results and they are not precise enough to estimate the accuracy of the amplitude. This is why we have performed a series of numerical tests and we have chosen to represent them from the value of $Q$ which is defined by:
$$Q=\frac{\max\nolimits_t\left|p(t,x,0)\right|_{full wave}}{\max\nolimits_t\left|p(t,x,0)\right|_{one-way}}.$$

\begin{figure}
\begin{center}
\includegraphics[scale=0.35]{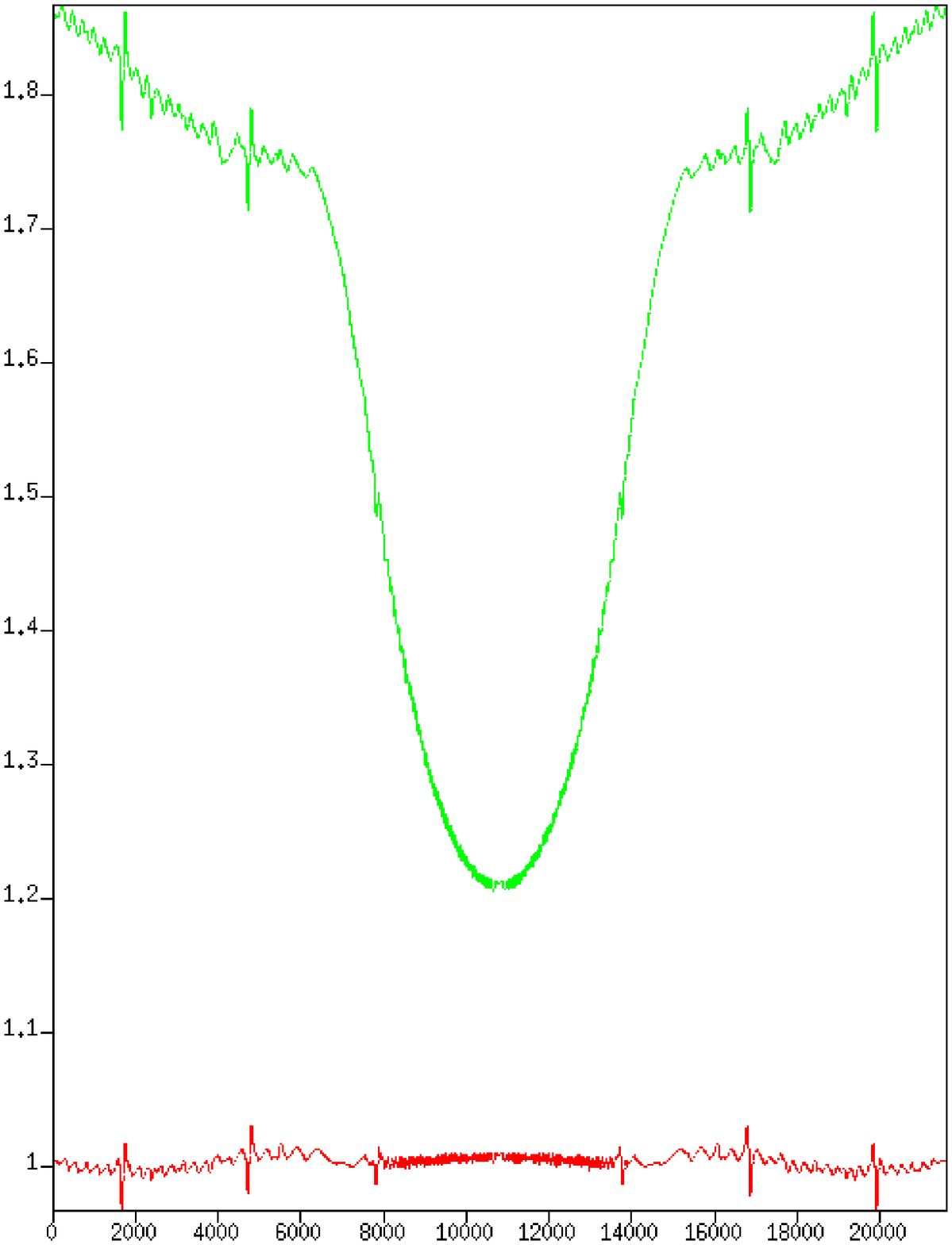}
\end{center}
\mcaption{Values of $Q(x)$ in VM 1. Top curve: one-way solution without transmission. Lower curve: one-way solution with transmission (rhs).}
\label{16002400_4vs0}
\end{figure}
On Fig. \ref{16002400_4vs0}, the lower curve depicts the value of $Q$ for the one-way solution with $\epsilon =0$. Its values are very close to $1$ which confirms what was observed on the corresponding seismogram. The top curve represents the values of $Q$ for the one-way solution without the transmission term. This picture strengthens the previous conclusion since the minimum value of the error is $20\%$. Hence it is essential to include the transmission into the model to get accurate amplitudes.\\
Now, being convinced that the transmission must be included, the theory shows that it can be equivalently introduced whether into the one-way equations or into the rhs of the system. The first way allows to write a system of one-way equations which are quite close to the equations considered by Zhang \emph{and al.} \cite{zhang} as an approximation of the full wave equation. The second one is the most natural in the formalism of Bremmer series and was introduced by Le Rousseau and De Hoop \cite{ref3,ref1}. \\
\newsavebox{\figtesta}
\begin{figure}
\begin{center}
\begin{pspicture}(15,8)
\savebox{\figtesta}(15,8)[h]{
\includegraphics[scale=0.35]{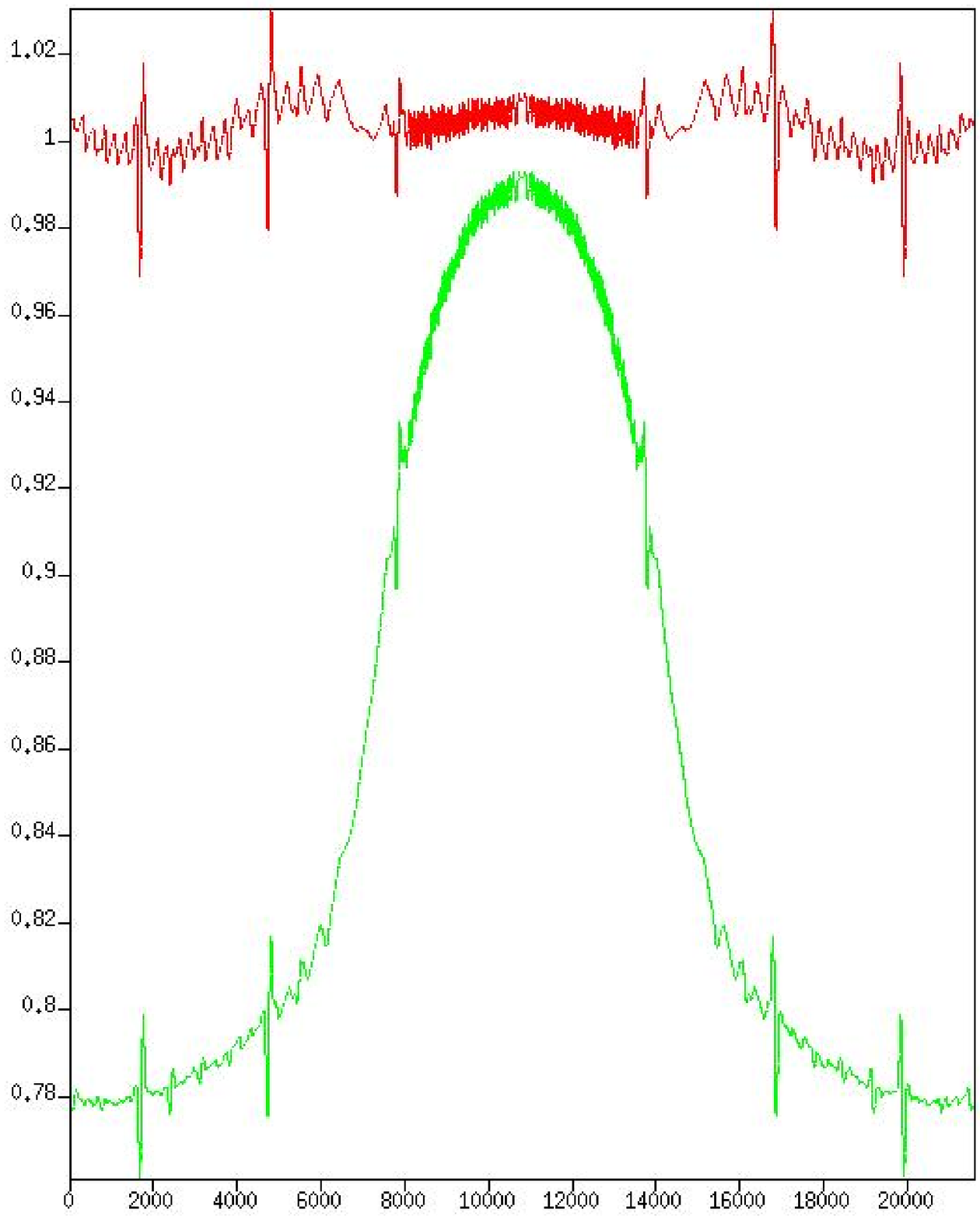}}
\rput[bl](0,0){\usebox{\figtesta}}
\end{pspicture}
\end{center}
\mcaption{Values of $Q(x)$ for the one-way solution with transmission in VM 1. Top curve: $\epsilon=0$. Lower curve: $\epsilon=1$}
\label{16002400}
\end{figure}
On Fig. \ref{16002400}, we have collected the values of $Q$ for $\epsilon=0$ and $\epsilon=1$.  %
When $\epsilon=1$, the results are correct in a neighbourhood of the shot and then they spoil quickly. Hence this numerical test indicates that for VM 1, the best approach consists in including the transmission into the rhs.\\
Thus the result for $ \epsilon=1$ may deteriorate more and more as fast as the velocity contrast is high. That motivates the next test on Fig. \ref{16005000} where we can observe that actually the results spoil more and more: the minimum value of the error is now $8\%$ versus $1\%$ for VM 1.\\

\newsavebox{\figtestb}
\begin{figure}
\begin{center}
\begin{pspicture}(15,8)
\savebox{\figtestb}(15,8)[h]{
\includegraphics[scale=0.35]{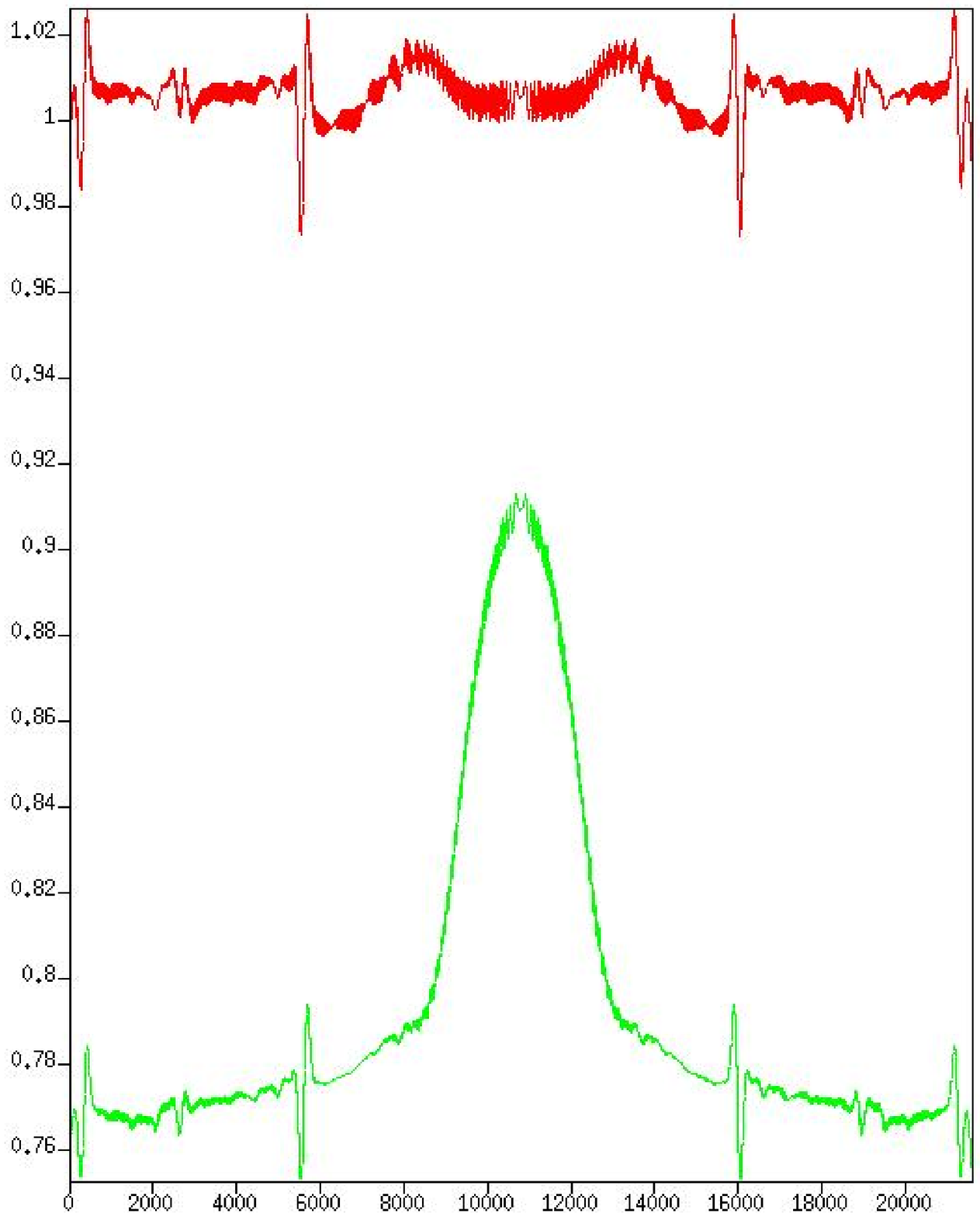}}
\rput[bl](0,0){\usebox{\figtestb}}
\end{pspicture}
\end{center}
\mcaption{Values of $Q(x)$ for the one-way solution with transmission in VM 2. Top curve: $\epsilon=0$. Lower curve: $\epsilon=1$}
\label{16005000}
\end{figure}
On Fig. \ref{5hhfchfac}, we have collected three curves obtained from the one-way solution with $\epsilon=1$. Each of them are distinguished from the value of the velocity contrast. The pattern shows that the lower the contrast is the larger the neigbourhood of the shot in which the results are satisfactory is. This figure seems to indicate that the model with $\epsilon=1$ is correct when the propagation medium is smooth.

\newsavebox{\figtestd}
\begin{figure}
\begin{center}
\begin{pspicture}(15,8)
\savebox{\figtestd}(15,8)[h]{
\includegraphics[scale=0.35]{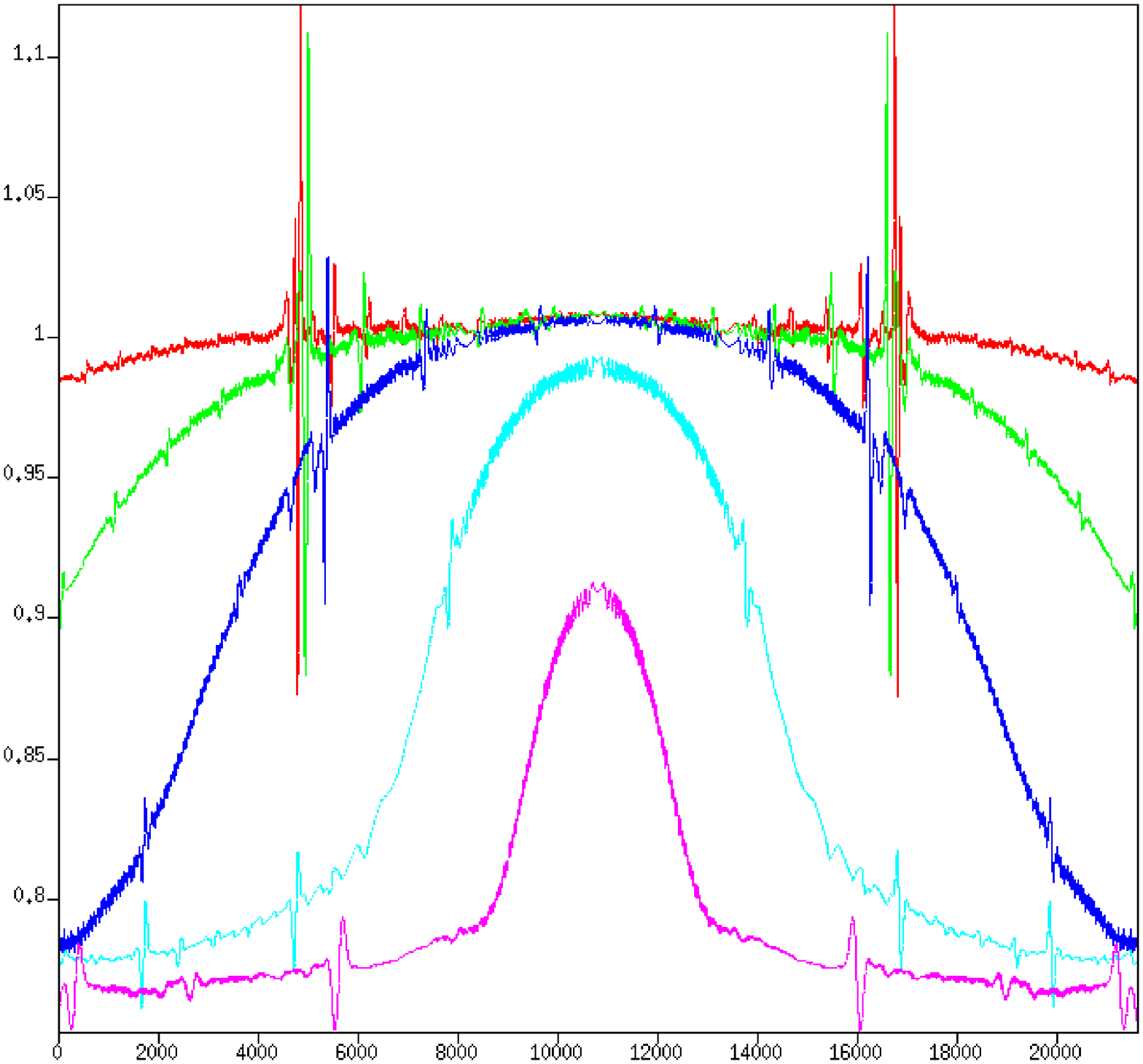}}
\rput[bl](0,0){\usebox{\figtestd}}
\uput{.2}[135](5.5,2.25){\tiny \colorbox{white}{1600/1800} \normalsize}
\uput{.2}[135](9.5,4){\tiny \colorbox{white}{1600/2400} \normalsize}
\uput{.2}[135](8,0.5){\tiny \colorbox{white}{1600/5000} \normalsize}
\uput{.2}[135](11.5,4.2){\tiny \colorbox{white}{1600/1700} \normalsize}
\uput{.2}[135](5.5,5){\tiny \colorbox{white}{1600/1650} \normalsize}
\end{pspicture}
\end{center}
\mcaption{Values of $Q(x)$ for the one-way solution with $\epsilon=1$ and with different velocity contrasts}
\label{5hhfchfac}
\end{figure}

\newsavebox{\figtestc}
\begin{figure}
\begin{center}
\begin{pspicture}(15,10)
\savebox{\figtestc}(15,10)[h]{
\includegraphics[scale=0.35]{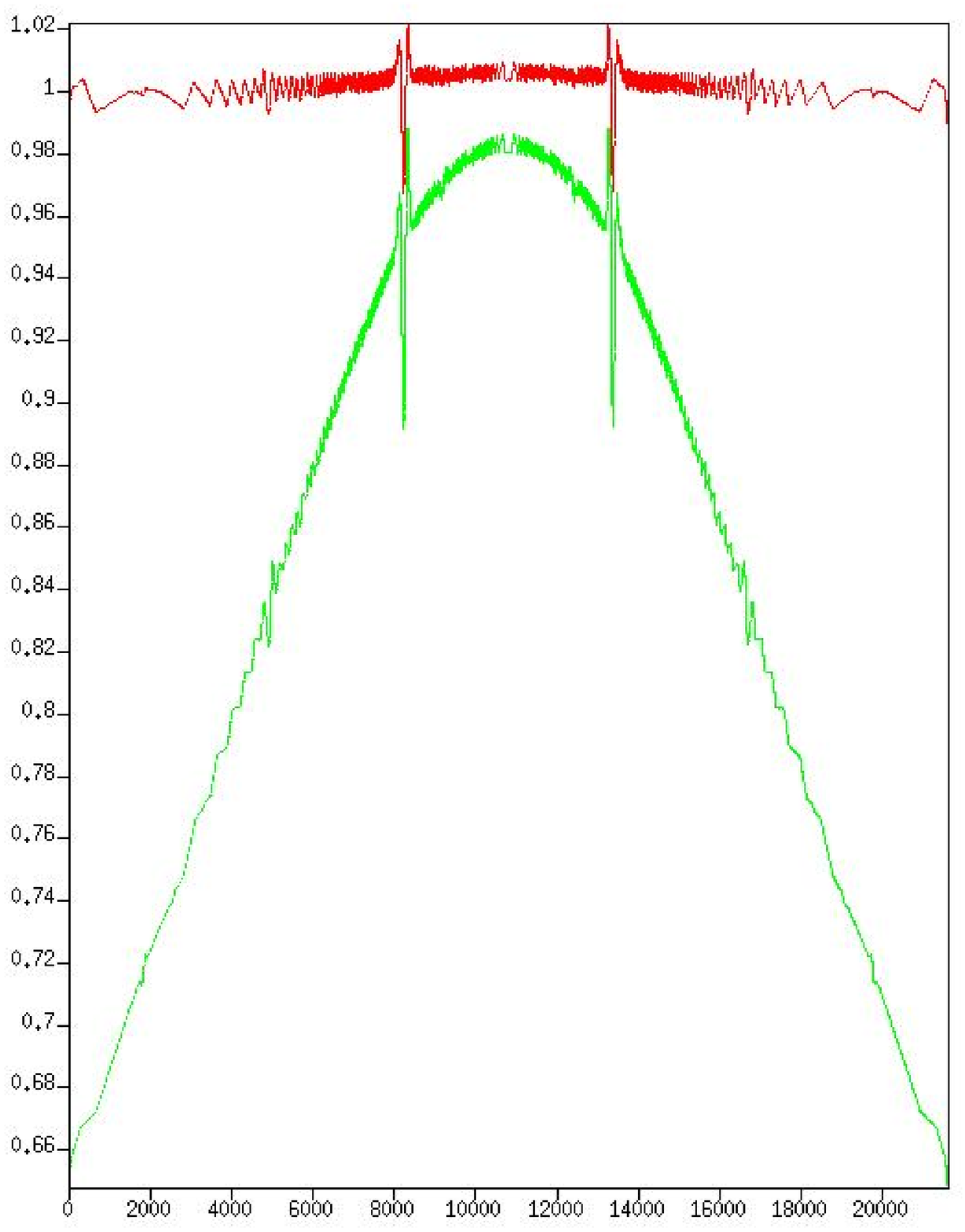}}
\rput[bl](0,0){\usebox{\figtestc}}
\end{pspicture}
\end{center}
\mcaption{Values of $Q(x)$ for the one-way solution with transmission in VM 3. Top curve: $\epsilon=0$. Lower curve: $\epsilon=1$}
\label{24001600}
\end{figure}
On Fig. \ref{24001600} we have depicted the results for VM 3. It is just to show that the same conclusion holds even if the upper velocity is larger than the lower one.
Let us mention that each curve shows some instabilities (depicted as local maxima) which are due to the periodicities created by the FFTs. 
\section{conclusion}
In this paper we consider the numerical analysis of two one-way systems derived from the general modeling of M. De Hoop \cite{ref1}.\\
Such a formulation is used to replace the full wave equation by a system of one-way wave equations whose computational burden is lower than the one associated to the finite difference solution of the wave equation. Moreover it permits to unpack and identify the multiples from the primary reflections. 
We have include a transmission term in the one-way model. Then numerical tests have been performed
in the 2D case and they indicate that accounting for the
transmission improves significantly the amplitude of the solution. 
The computational algorithm has been optimized in such a way that
its complexity is now of the same order than the one of the
one-way solver without transmission. Hence to add the correcting
term does not penalize the computational method. Now the solution
is based on a method of propagators whose 
numerical approximation accuracy is very sensitive to the position of this
correcting term. To put it into the one-way system, which is the
case $\varepsilon=1$, seems to deteriorate the results. But this
approach should be interesting since it corresponds to the same
idea than the one proposed by Zhang \emph{et al.} in \cite{zhang} where
time-arrivals are computed from the solution of a second-order
wave equation obtained by factoring the full wave equation. 

In this paper, we have performed a numerical test which shows that the accuracy of the method with $\epsilon=1$ is similar to the one of the method with $\epsilon=0$ when the velocity contrasts decrease.\\
In the proposed methods, the transmission operator is approximated by a zero order pseudodifferential operator which is exact for stratified media. A higher order operator, that should account for media with lateral velocity variations, is the topic of a current research \cite{badupra2}. 
\clearpage

\bibliographystyle{plain}
\bibliography{biblio}

\end{document}